\documentclass{article}
\usepackage[top=30truemm, bottom=30truemm, left=25truemm, right=25truemm]{geometry}

\usepackage{amsmath, amssymb, amsthm}

\newtheorem{defi}{Definition}[section]
\newtheorem{nota}[defi]{Notation}
\newtheorem{lem}[defi]{Lemma}
\newtheorem{prop}[defi]{Proposition}
\newtheorem{thm}[defi]{Theorem}
\newtheorem{cor}[defi]{Corollary}

\allowdisplaybreaks

\title{Terwilliger algebras of generalized wreath products of association schemes}
\author{Yuta Watanabe}
\date{September 2024}

\begin{document}

\maketitle

\begin{abstract}
The generalized wreath product of symmetric association schemes was introduced by R.~A.~Bailey in the European Journal of Combinatorics 27 (2006) 428--435.
It is recognized as a unification of both the wreath product and the direct product of symmetric association schemes. While its potential applicability to any association scheme had been implied, this paper provides a formal and explicit confirmation of that claim.
Moreover, we establish the irreducible representations of its adjacency algebra and Terwilliger algebra.

\end{abstract}

\section{Introduction}

Association schemes were initially introduced in statistics within the framework of experimental design and can be regarded as an extension of groups from a combinatorial standpoint.
Specifically, an association scheme naturally emerges when a group acts transitively on a finite set.
In recent years, building on the foundational work of Delsarte \cite{D}, these schemes have found extensive applications in various combinatorial domains, including coding theory and design theory.
In 1992, the Terwilliger algebra was first introduced in \cite{Ter} as a novel approach to studying commutative association schemes.
The Terwilliger algebra has been examined across diverse classes of association schemes.
For example, investigations have focused on group association schemes \cite{BM}, direct products of association schemes \cite{HKM}, wreath products of one-class association schemes \cite{BST}, wreath products of association schemes with some assumptions \cite{HKM}, wreath products of quasi-thin schemes \cite{K}, and wreath products of arbitrary association schemes \cite{MX}, among others.
For further information, we direct readers to \cite{Bai04, BBIT, Z} and the references contained therein.
The generalized wreath product of symmetric association schemes was introduced in \cite{Bai}.
It functions as a combinatorial counterpart to the generalized wreath product of transitive permutation groups, which was introduced in \cite{BPRS}.
Despite its designation, the generalized wreath product extends both the wreath product and the direct product, with each product manifesting as a special case.

This paper generalizes the concept of the ``generalized wreath product'' to arbitrary association schemes. 
Since the work of \cite{Bai} did not critically depend on symmetry or commutativity, they implicitly suggested that the generalized wreath product could be defined for any association scheme.
This paper provides a formal and explicit affirmation of this notion (see Theorem \ref{thm:X}).
Furthermore, it derives the irreducible representations of the adjacency algebra, a task previously undertaken by \cite{Bai} in the symmetric case (see Theorems \ref{thm:A} and \ref{thm:E}). 
As the proofs in \cite{Bai} heavily relied on commutativity, the results presented here significantly extend the scope of previous research.
Additionally, this paper determines the irreducible representations of the Terwilliger algebra, an issue that remains unresolved even in the symmetric case (see Theorems \ref{thm:J}, \ref{thm:eps}, \ref{thm:T} and \ref{thm:epsilon}).
The existing results concerning the Terwilliger algebras of wreath products and direct products from earlier works \cite{BST, HKM, K, MX} are shown to be special instances of our findings.

The structure of this paper is organized as follows.
In Section 2, we review the fundamental notation of association schemes.
In Section 3, we introduce the generalized wreath products for arbitrary association schemes.
In Section 4, we express the adjacency algebra in terms of those of its constituent factors.
In Section 5, we identify the central primitive idempotents of the adjacency algebra.
In Section 6, we represent the dual adjacency algebra in terms of the corresponding factors.
In Section 7, we formulate the Terwilliger algebra based on the structure of each factor.
In Section 8, we express the principal two-sided ideal generated by the all-ones matrix and determine the primitive idempotents of the Terwilliger algebra.
In Section 9, we establish the central primitive idempotents of the Terwilliger algebra.
In Section 10, we discuss the application of main theorems to the special cases studied in \cite{HKM, MX}.

We conclude this section by introducing some notation pertaining to matrix algebras, which will be utilized throughout this paper.
For a non-empty finite set $X$, let $\operatorname{Mat}_X(\mathbb{C})$ denote the algebra over the complex field $\mathbb{C}$ consisting of matrices with rows and columns indexed by $X$, with entries in $\mathbb{C}$.
Let $I_X \in \operatorname{Mat}_X(\mathbb{C})$ denote the identity matrix and let $J_X \in \operatorname{Mat}_X(\mathbb{C})$ denote the all-ones matrix.
For a matrix $U \in \operatorname{Mat}_X(\mathbb{C})$, let $\langle U \rangle$ denote the subspace of $\operatorname{Mat}_X(\mathbb{C})$ spanned by $U$.
For a matrix $U \in \operatorname{Mat}_X(\mathbb{C})$ and a subspace $\mathcal{A}$ of $\operatorname{Mat}_X(\mathbb{C})$,
define $U\mathcal{A} = \{UA \mid A \in \mathcal{A}\}$ and $\mathcal{A}U = \{AU \mid A \in \mathcal{A}\}$.
Both $U\mathcal{A}$ and $\mathcal{A}U$ are subspaces of $\operatorname{Mat}_X(\mathbb{C})$.
Let $Y$ be another non-empty finite set.
The Kronecker product of matrices $U \in \operatorname{Mat}_X(\mathbb{C})$ and $W \in \operatorname{Mat}_Y(\mathbb{C})$ is denoted by $U \otimes W \in \operatorname{Mat}_{X \times Y}(\mathbb{C})$.
For subspaces $\mathcal{A}$ of $\operatorname{Mat}_X(\mathbb{C})$ and $\mathcal{B}$ of $\operatorname{Mat}_Y(\mathbb{C})$,
$\mathcal{A} \otimes \mathcal{B} = \operatorname{Span}_{\mathbb{C}}\{A \otimes B \mid A \in \mathcal{A}, B \in \mathcal{B}\}$ is also a subspace of $\operatorname{Mat}_{X \times Y}(\mathbb{C})$.
For a matrix $U \in \operatorname{Mat}_X(\mathbb{C})$ and a subspace $\mathcal{A}$ of $\operatorname{Mat}_Y(\mathbb{C})$,
define $U \otimes \mathcal{A} = \langle U \rangle \otimes \mathcal{A}$ and $\mathcal{A} \otimes U = \mathcal{A} \otimes \langle U \rangle$.
Then, $U \otimes \mathcal{A}$ and $\mathcal{A} \otimes U$ are subspaces of $\operatorname{Mat}_{X \times Y}(\mathbb{C})$ and $\operatorname{Mat}_{Y \times X}(\mathbb{C})$, respectively.

\section{Association schemes}

In this section, we provide a concise review of the notation and fundamental concepts related to association schemes.
For further details, the reader is directed to \cite{Bai04, BBIT, Z} and the references cited therein.

Let $X$ be a finite set with at least two elements.
A pair $\mathfrak{X} = (X,\{R_i\}_{i=0}^d)$, consisting of $X$ and a collection of subsets of $X \times X$, is referred to as a $d$-class association scheme if it satisfies the following four conditions:
\begin{enumerate}
\item[(R1)] $R_0 = \{(x,x) \mid x \in X\}$;
\item[(R2)] $\{R_i\}_{i=0}^d$ forms a partition of $X \times X$;
\item[(R3)] For $0 \le i \le d$, there exists $0 \le i' \le d$ such that $R_{i'} = \{(y,x) \mid (x,y) \in R_i\}$;
\item[(R4)] For $0 \le i,j,k \le d$, and a given $(x,y) \in R_k$, the value $\mathsf{p}_{i,j}^k = |\{z \in X \mid (x,z) \in R_i, (z,y) \in R_j\}|$ depends only on $i,j,k$. The number $\mathsf{p}_{i,j}^k$ is called the intersection number.
\end{enumerate}
Additionally, if the following condition holds, a pair $\mathfrak{X} = (X,\{R_i\}_{i=0}^d)$ is termed a commutative association scheme.
\begin{enumerate}
\item[(R5)] For $0 \le i,j,k \le d$, the numbers in condition (R4) satisfy $\mathsf{p}_{i,j}^k = \mathsf{p}_{j,i}^k$.
\end{enumerate}
Moreover, if the following condition holds, a pair $\mathfrak{X} = (X,\{R_i\}_{i=0}^d)$ is known as a symmetric association scheme.
\begin{enumerate}
\item[(R6)] For $0 \le i \le d$, we have $i=i'$ in the condition (R3).
\end{enumerate}
We leave it to the reader to verify that every symmetric association scheme is commutative.
We define the following:
\begin{align*}
\mathsf{k}_i = \mathsf{p}_{i,i'}^0
&&
0 \le i \le d.
\end{align*}
Note that $\mathsf{k}_i$ denotes the number of $y \in X$ such that $(x,y) \in R_i$ for some $x \in X$.
We refer to $\{\mathsf{k}_i\}_{i=0}^d$ as the valencies of $\mathfrak{X}$. 

Let $\mathfrak{X} = (X,\{R_i\}_{i=0}^d)$ be a $d$-class association scheme.
For $0 \le i \le d$, define a matrix $A_i \in \operatorname{Mat}_X(\mathbb{C})$ by
\begin{align*}
(A_i)_{x,y} =
\begin{cases}
1 & \text{if $(x,y) \in R_i$}, \\ 0 & \text{if $(x,y) \not\in R_i$}
\end{cases}
&&
x,y \in X.
\end{align*}
We refer to $\{A_i\}_{i=0}^d$ as the adjacency matrices of $\mathfrak{X}$.
The subspace $\mathcal{A} \subset \operatorname{Mat}_X(\mathbb{C})$ spanned by $\{A_i\}_{i=0}^d$, by conditions (R1),(R3),(R4), forms a semi-simple subalgebra, known as the adjacency algebra of $\mathfrak{X}$.
Moreover, $\{A_i\}_{i=0}^d$ provides a basis for $\mathcal{A}$.
Note that $J_X \in \mathcal{A}$ by condition (R2) and $\langle J_X \rangle$ constitutes a two-sided ideal of $\mathcal{A}$.
As the adjacency algebra $\mathcal{A}$ is semi-simple, it can be decomposed as a direct sum of central simple algebras:
\[
\mathcal{A} = \bigoplus_{i=0}^{l} \mathcal{A}_i.
\]
Without loss of generality, assume $\mathcal{A}_0 = \langle J_X \rangle$.
The identity element of $\mathcal{A}$ is expressed as
\begin{align*}
I_X = E_0 + E_1 + \cdots + E_l,
&&
E_j \in \mathcal{A}_j.
\end{align*}
Since $\mathcal{A}_0 = \langle J_X \rangle$, we have
\[
E_0 = \frac{1}{|X|}J_X.
\]
We refer to $\{E_j\}_{j=0}^l$ as the primitive central idempotents of $\mathcal{A}$.

Fix $x \in X$.
For $0 \le i \le d$, let $E_i^* = E_i^*(x)$ denote the diagonal matrix in $\operatorname{Mat}_X(\mathbb{C})$ defined by
\begin{align*}
(E_i^*)_{y,y} = (A_i)_{x,y} && y \in X.
\end{align*}
We refer to $\{E_i^*\}_{i=0}^d$ as the dual primitive idempotents of $\mathfrak{X}$.
The subspace $\mathcal{A}^* = \mathcal{A}^*(x)$ of $\operatorname{Mat}_X(\mathbb{C})$ spanned by $\{E_i^*\}_{i=0}^d$ forms a commutative semi-simple subalgebra, known as the dual adjacency algebra of $\mathfrak{X}$.
The Terwilliger algebra $\mathcal{T} = \mathcal{T}(x)$ of $\mathfrak{X}$ is the subalgebra of $\operatorname{Mat}_X(\mathbb{C})$ generated by
both the adjacency algebra $\mathcal{A}$ and the dual adjacency algebra $\mathcal{A}^*$.

\begin{lem}\label{lem}
With the above notation, the following statemanets hold:
\begin{enumerate}
\item $J_X\mathcal{T} = J_X\mathcal{A}^*$ and $\mathcal{T}J_X = \mathcal{A}^*J_X$.
\item $E_0^*\mathcal{T} = E_0^*\mathcal{A}$ and $\mathcal{T}E_0^* = \mathcal{A}E_0^*$.
\end{enumerate}
\end{lem}
\begin{proof}
\begin{enumerate}
\item The inclusion $\mathcal{A}^* \subset \mathcal{T}$ implies $J_X\mathcal{A}^* \subset J_X\mathcal{T}$.
Since $\mathcal{A}_0 = \langle J_X \rangle$ is a two-sided ideal of $\mathcal{A}$ and $I_X \in \mathcal{A}^*$, it follows that $J_X\mathcal{A} = \mathcal{A}_0 \subset J_X\mathcal{A}^*$.
Thus, the inclusions show that $J_X\mathcal{T} = J_X\mathcal{A}^*$.
Taking the transpose of both sides of this equality yields $\mathcal{T}J_X = \mathcal{A}^*J_X$.
\item The inclusion $\mathcal{A} \subset \mathcal{T}$ implies $E_0^*\mathcal{A} \subset E_0^*\mathcal{T}$.
Since $\langle E_0^* \rangle$ is a two-sided ideal of $\mathcal{A}^*$ and $I_X \in \mathcal{A}$, it follows that $E_0^*\mathcal{A}^* \subset E_0^*\mathcal{A}$.
Thus, the inclusions show that $E_0^*\mathcal{T} = E_0^*\mathcal{A}$.
Taking the transpose of both sides of this equality yields $\mathcal{T}E_0^* = \mathcal{A}E_0^*$.
\end{enumerate}
\end{proof}

Let $\mathcal{J} = \mathcal{J}(x)$ denote the principal two-sided ideal of $\mathcal{T}$ generated by the all-ones matrix $J_X$.
In other words, we define $\mathcal{J} = \mathcal{T}J_X\mathcal{T}$.
By Lemma \ref{lem}, it can be expressed as $\mathcal{J} = \mathcal{A}^*J_X\mathcal{A}^*$.
Since the Terwilliger algebra $\mathcal{T}$ is semi-simple, it can be decomposed into a direct sum of central simple algebras:
\[
\mathcal{T} = \bigoplus_{i=0}^{t} \mathcal{T}_i.
\]
Without loss of generality, assume $\mathcal{T}_0 = \mathcal{J}$.
Let the identity element of $\mathcal{T}$ be expressed as
\begin{align*}
I_X = \varepsilon_0 + \varepsilon_1 + \cdots + \varepsilon_t,
&&
\varepsilon_k \in \mathcal{T}_k.
\end{align*}
We refer to $\{\varepsilon_k\}_{k=0}^t$ the primitive central idempotents of $\mathcal{T}$.
Since we assume $\mathcal{T}_0 = \mathcal{J}$, it follows that
\[
\varepsilon_0 = \sum_{i=0}^d \frac{1}{\mathsf{k}_i} E_i^* J_X E_i^*.
\]
We term this as the principal idempotent of $\mathcal{T}$.

\section{Generalized wreath products of association schemes}

In this section, we succinctly revisit the notation of partially ordered sets (posets) and present generalized wreath products of association schemes. 
For additional details on posets, we direct the reader to \cite{Sta}.

Let $P$ be a non-empty finite poset.
An anti-chain in $P$ refers to a subset of $P$ in which no two distinct elements are comparable.
The empty set $\emptyset$ is also considered an anti-chain.
A down-set in $P$ is a subset $D$ of $P$ such that if $p \in D$ and $q < p$, then $q \in D$.
Similarly, an up-set in $P$ is a subset $U$ of $P$ such that if $p \in U$ and $p < q$, then $q \in U$.
For an anti-chain $M$ in $P$, we define
\begin{align*}
M^\downarrow &= \{ p \in P \mid \text{there exists $q \in M$ such that $p < q$}\},\\ M^\uparrow &= \{ p \in P \mid \text{there exists $q \in M$ such that $q < p$}\}. \end{align*}
Then $M^\downarrow$ constitutes a down-set in $P$, and $M^\uparrow$ constitutes an up-set in $P$.
A subset $C$ of $P$ is termed convex if $q \in C$ whenever $p < q < r$ with $p, r \in C$. Note that anti-chains, down-sets, and up-sets are all convex.
For a subset $S$ of $P$, let $\max(S)$ (resp.\ $\min(S)$) denote the set of all maximal (resp.\ minimal) elements in $S$.
It is noteworthy that $\max(S)$ and $\min(S)$ are anti-chains in $P$ for any subset $S$ of $P$.

For each $p \in P$, consider a finite set $X_p$ with at least two elements and a $d_p$-class association scheme $\mathfrak{X}_p = (X_p, \{R_{p,i}\}_{i=0}^{d_p})$.
Recall condition (R3), which states that for each $0 \le i \le d_p$, there exists $0 \le i' \le d_p$ such that $R_{p,i'} = \{(y,x) \mid (x,y) \in R_{p,i}\}$.
Let $\{\mathsf{p}_{p,i,j}^k\}_{i,j,k=0}^{d_p}$ and $\{\mathsf{k}_{p,i}\}_{i=0}^{d_p}$ denote the intersection numbers and the valencies for $\mathfrak{X}_p$, respectively.
We consider the Cartesian product $X = \prod_{p \in P} X_p$.
For $x = (x_p)_{p \in P}$ and $y = (y_p)_{p \in P}$ in $X$, define
\[
M = \max\{p \in P \mid x_p \neq y_p\}.
\]
Since $M$ represents a set of maximal elements, it constitutes an anti-chain in $P$.
For each $p \in M$, by conditions (R1) and (R2), there uniquely exists $1 \le i_p \le d_p$ such that $(x_p, y_p) \in R_{p,i_p}$.
In this case, we say that $x$ is $i$-adjacent to $y$, where $i = (i_p)_{p \in M}$.
Let $\mathbb{I}$ denote the set of all such $i$. In other words,
\[
\mathbb{I} = \{ (i_p)_{p \in M} \mid \text{$M$ is an anti-chain in $P$}, 1 \le i_p \le d_p \; (p \in M)\}.
\]

\begin{lem}\label{lem:i'}
Let $i \in \mathbb{I}$ and $x,y \in X$.
Write $i = (i_p)_{p \in M}$ and define $i' = (i'_p)_{p \in M}$.
Then $x$ is $i$-adjacent to $y$ if and only if $y$ is $i'$-adjacent to $x$.
\end{lem}
\begin{proof}
Straightforward.
\end{proof}

\begin{lem}\label{lem:i}
Let $i \in \mathbb{I}$ and $x, y \in X$.
Write $i = (i_p)_{p \in M}$ and $x = (x_p)_{p \in P}$, $y = (y_p)_{p \in P}$.
Then $x$ is $i$-adjacent to $y$ if and only if for $p \in P$, the following conditions are satisfied:
\begin{enumerate}
\item If $p \in M$, then $(x_p, y_p) \in R_{p,i_p}$. In particular, $x_p \neq y_p$.
\item If $p \not\in M \cup M^\downarrow$, then $x_p = y_p$.
\end{enumerate}
\end{lem}
\begin{proof}
Assume $x$ is $i$-adjacent to $y$.
Condition (i) follows directly from the definition.
To verify condition (ii), assume $x_p \neq y_p$ and demonstrate that $p \in M \cup M^\downarrow$.
By the definition of $M$, there exists $q \in M$ such that $p \le q$, which implies $p \in M \cup M^\downarrow$.

Conversely, assume conditions (i) and (ii) are satisfied.
Define $Q = \{ p \in P \mid x_p \neq y_p \}$.
Given conditions (i) and (ii), it follows that $M \subset Q \subset M \cup M^\downarrow$, which implies $M = \max(Q)$.
Therefore, by assumption (i), $x$ is $i$-adjacent to $y$.
\end{proof}

\begin{lem}\label{lem:pijk1}
Let $i,j,k \in \mathbb{I}$ and $x,y \in X$ such that $x$ is $k$-adjacent to $y$.
Define $i = (i_p)_{p \in M}$, $j = (j_p)_{p \in N}$, and $k = (k_p)_{p \in L}$, and assume the following conditions hold:
\begin{enumerate}
\item $L \cup L^\downarrow \subset (M \cup M^\downarrow) \cup (N \cup N^\downarrow)$,
\item $M \setminus (N \cup N^\downarrow) \subset L$,
\item $N \setminus (M \cup M^\downarrow) \subset L$,
\item $i_p = j'_p$ if $p \in (M \cap N) \setminus (L \cup L^\downarrow)$,
\item $i_p = k_p$ if $p \in (M \cap L) \setminus (N \cup N^\downarrow)$,
\item $j_p = k_p$ if $p \in (N \cap L) \setminus (M \cup M^\downarrow)$.
\end{enumerate}
Then the number of $z \in X$ such that $x$ is $i$-adjacent to $z$ and $z$ is $j$-adjacent to $y$ is given by
\[
\prod_{p \in M^\downarrow \cap N^\downarrow} |X_p| \times
\prod_{p \in M^\downarrow \cap N} \mathsf{k}_{p,j_p} \times
\prod_{p \in M \cap N^\downarrow} \mathsf{k}_{p,i_p} \times
\prod_{p \in M \cap N \cap L} \mathsf{p}_{p,i_p,j_p}^{k_p} \times 
\prod_{p \in (M \cap N) \setminus (L \cup L^\downarrow)} \mathsf{k}_{p,i_p}.
\]
\end{lem}

\begin{proof}

From Lemma \ref{lem:i}, the adjacency conditions for two points $x$ and $y$ can be independently expressed in terms of $x_p$ and $y_p$ for each $p$.
Consequently, it is sufficient to count the suitable $z_p$ for each $p$ and compute their product.
Specifically, let $\mathsf{n}_p$ denote the number of $z_p$ satisfying the following conditions:
(1) If $p \in M$, then $(x_p, z_p) \in R_{p, i_p}$;
(2) If $p \not\in M \cup M^\downarrow$, then $x_p = z_p$.
(3) If $p \in N$, then $(z_p, y_p) \in R_{p, j_p}$;
(4) If $p \not\in N \cup N^\downarrow$, then $z_p = y_p$. 
Thus, the product $\prod_{p \in P} \mathsf{n}_p$ represents the quantity we seek.
According to the relationship between $p$ and $M$, the case is categorized into the following three scenarios: (a) $p \in M^\downarrow$; (b) $p \in M$; (c) $p \not\in M \cup M^\downarrow$.
Furthermore, based on the relationship between $p$ and $N$, each scenario is further divided into three sub-cases: ($\alpha$) $p \in N^\downarrow$; ($\beta$) $p \in N$; ($\gamma$) $p \not\in N \cup N^\downarrow$.
Thus, a total of nine distinct patterns are considered.
\begin{enumerate}
\item[(a-$\alpha$)]
Consider the case where $p \in M^\downarrow$ and $p \in N^\downarrow$.
Here, $\mathsf{n}_p$ represents the number of $z_p$ without any additional conditions.
Thus, $\mathsf{n}_p = |X_p|$.
\item[(a-$\beta$)]
Consider the case where $p \in M^\downarrow$ and $p \in N$.
Here, $\mathsf{n}_p$ represents the number of $z_p$ that satisfy $(z_p, y_p) \in R_{p, j_p}$.
Therefore, $\mathsf{n}_p = \mathsf{k}_{p, j_p'} = \mathsf{k}_{p, j_p}$.
\item[(a-$\gamma$)]
Consider the case where $p \in M^\downarrow$ and $p \not\in N \cup N^\downarrow$.
Here, $\mathsf{n}_p$ represents the number of $z_p$ that satisfy $z_p = y_p$.
Thus, $\mathsf{n}_p = 1$.
\item[(b-$\alpha$)]
Consider the case where $p \in M$ and $p \in N^\downarrow$.
Here, $\mathsf{n}_p$ represents the number of $z_p$ that satisfy $(x_p,z_p) \in R_{p,i_p}$.
Therefore, $n_p = \mathsf{k}_{p,i_p}$.
\item[(b-$\beta$)]
Consider the case where $p \in M$ and $p \in N$.
Here, $\mathsf{n}_p$ represents the number of $z_p$ satisfying both $(x_p, z_p) \in R_{p, i_p}$ and $(z_p, y_p) \in R_{p, j_p}$. According to assumption (i), either $p \in L$ or $p \not\in L \cup L^\downarrow$.
If $p \in L$, then by Lemma \ref{lem:i}, $(x_p, y_p) \in R_{p, k_p}$.
Thus, $\mathsf{n}_p = \mathsf{p}_{p, i_p, j_p}^{k_p}$.
If $p \not\in L \cup L^\downarrow$, then by Lemma \ref{lem:i}, $x_p = y_p$. 
Hence, $\mathsf{n}_p = \mathsf{p}_{p, i_p, j_p}^{0} = \mathsf{k}_{p, i_p}$ by assumption (iv).
\item[(b-$\gamma$)]
Consider the case where $p \in M$ and $p \not\in N \cup N^\downarrow$.
Here, $\mathsf{n}_p$ represents the number of $z_p$ satisfying $(x_p, z_p) \in R_{p, i_p}$ and $z_p = y_p$.
According to assumptions (ii) and (v), $p \in L$ and $i_p = k_p$.
By Lemma \ref{lem:i}, this implies $(x_p, y_p) \in R_{p, k_p}$.
Thus, $\mathsf{n}_p = 1$.
\item[(c-$\alpha$)]
Consider the case where $p \not\in M \cup M^\downarrow$ and $p \in N^\downarrow$.
Here, $\mathsf{n}_p$ represents the number of $z_p$ that satisfy $x_p = z_p$.
Thus, $\mathsf{n}_p = 1$.
\item[(c-$\beta$)]
Consider the case where $p \not\in M \cup M^\downarrow$ and $p \in N$.
Here, $\mathsf{n}_p$ represents the number of $z_p$ that satisfy $x_p = z_p$ and $(z_p,y_p) \in R_{p,j_p}$.
According to assumptions (iii) and (vi), $p \in L$ and $j_p = k_p$.
By Lemma \ref{lem:i}, this implies $(x_p, y_p) \in R_{p,k_p}$.
Thus, $\mathsf{n}_p = 1$.
\item[(c-$\gamma$)]
Consider the case where $p \not\in M \cup M^\downarrow$ and $p \not\in N \cup N^\downarrow$.
Here, $\mathsf{n}_p$ represents the number of $z_p$ that satisfy $x_p = z_p$ and $z_p = y_p$.
According to assumption (i), $p \not\in L \cup L^\downarrow$.
By Lemma \ref{lem:i}, this implies $x_p = y_p$.
Thus, $\mathsf{n}_p = 1$.
\end{enumerate}
The product of $\mathsf{n}_p$ across all cases yields the desired value.
\end{proof}

\begin{lem}\label{lem:pijk2}
Let $i, j, k \in \mathbb{I}$ and $x, y \in X$ such that $x$ is $k$-adjacent to $y$.
Suppose that at least one of the conditions (i)--(vi) in Lemma \ref{lem:pijk1} fails to hold.
Then, there does not exist a $z \in X$ such that $x$ is $i$-adjacent to $z$ and $z$ is $j$-adjacent to $y$.
\end{lem}
\begin{proof}

For each $p \in P$, let $\mathsf{n}_p$ denote the number of $z_p$ satisfying the following conditions:
(1) If $p \in M$, then $(x_p, z_p) \in R_{p, i_p}$;
(2) If $p \not\in M \cup M^\downarrow$, then $x_p = z_p$.
(3) If $p \in N$, then $(z_p, y_p) \in R_{p, j_p}$;
(4) If $p \not\in N \cup N^\downarrow$, then $z_p = y_p$.
As demonstrated in the proof of Lemma \ref{lem:pijk1}, it suffices to show that there exists $p \in P$ such that $\mathsf{n}_p = 0$. Let $x = (x_p)_{p \in P}$ and $y = (y_p)_{p \in P}$.
\begin{enumerate}
\item
Assume that condition (i) in Lemma \ref{lem:pijk1} does not hold.
We have either $L \not\subset (M \cup M^\downarrow) \cup (N \cup N^\downarrow)$ or $L^\downarrow \not\subset (M \cup M^\downarrow) \cup (N \cup N^\downarrow)$.
If $L \not\subset (M \cup M^\downarrow) \cup (N \cup N^\downarrow)$, there exists $p \in L$ such that $p \not\in M \cup M^\downarrow$ and $p \not\in N \cup N^\downarrow$.
Here, $\mathsf{n}_p$ represents the number of $z_p$ satisfying $z_p = x_p$ and $z_p = y_p$
Hence, $\mathsf{n}_p = 0$ because $x_p \neq y_p$ according to Lemma \ref{lem:i}.
If $L^\downarrow \not\subset (M \cup M^\downarrow) \cup (N \cup N^\downarrow)$, there exists $p \in L^\downarrow$ such that $p \not\in M \cup M^\downarrow$ and $p \not\in N \cup N^\downarrow$.
Consequently, there exists $q \in L$ such that $p < q$.
This $q$ satisfies $q \not\in M \cup M^\downarrow$ and $q \not\in N \cup N^\downarrow$.
As in the previous case, we have $\mathsf{n}_p = 0$.
\item
Assume that condition (ii) in Lemma \ref{lem:pijk1} does not hold.
There exists $p \in M$ such that $p \not\in N \cup N^\downarrow$ and $p \not\in L$.
We then have either $p \in L^\downarrow$ or $p \not\in (L \cup L^\downarrow)$.
If $p \in L^\downarrow$, there exists $q \in L$ such that $p < q$.
This $q$ satisfies $q \not\in M \cup M^\downarrow$ and $q \not\in N \cup N^\downarrow$.
As in case (i), we have $\mathsf{n}_p = 0$.
If $p \not\in (L \cup L^\downarrow)$, by Lemma \ref{lem:i}, $x_p = y_p$.
Here, $\mathsf{n}_p$ represents the number of $z_p$ satisfying $(x_p, z_p) \in R_{p,i_p}$ and $z_p = y_p$.
Hence, $\mathsf{n}_p = 0$ because $i_p \neq 0$.
\item
Assume that condition (iii) in Lemma \ref{lem:pijk1} does not hold.
There exists $p \in N$ such that $p \not\in M \cup M^\downarrow$ and $p \not\in L$.
We then have either $p \in L^\downarrow$ or $p \not\in (L \cup L^\downarrow)$.
If $p \in L^\downarrow$, there exists $q \in L$ such that $p < q$.
This $q$ satisfies $q \not\in M \cup M^\downarrow$ and $q \not\in N \cup N^\downarrow$.
As in case (i), we have $\mathsf{n}_p = 0$.
If $p \not\in (L \cup L^\downarrow)$, by Lemma \ref{lem:i}, $x_p = y_p$.
Here, $\mathsf{n}_p$ represents the number of $z_p$ satisfying $x_p = z_p$ and $(z_p, y_p) \in R_{p,j_p}$.
Hence, $\mathsf{n}_p = 0$ because $j_p \neq 0$.
\item 
Assume the condition (iv) in Lemma \ref{lem:pijk1} does not hold.
There exists $p \in (M \cap N) \setminus (L \cup L^\downarrow)$ such that $i_p \neq j_p'$.
Here, $\mathsf{n}_p$ represents the number of $z_p$ satisfying $(x_p, z_p) \in R_{p,i_p}$ and $(z_p, y_p) \in R_{p,j_p}$.
Hence, $\mathsf{n}_p = 0$ because $x_p = y_p$ according to Lemma \ref{lem:i}.
\item
Assume the condition (v) in Lemma \ref{lem:pijk1} does not hold.
There exists $p \in (N \cap L) \setminus (M \cup M^\downarrow)$ such that $j_p \neq k_p$.
Here, $\mathsf{n}_p$ represents the number of $z_p$ satisfying $z_p = x_p$ and $(z_p, y_p) \in R_{p,j_p}$.
Hence, $\mathsf{n}_p = 0$ because $(x_p, y_p) \in R_{p,k_p}$ according to Lemma \ref{lem:i}.
\item
Assume the condition (vi) in Lemma \ref{lem:pijk1} does not hold.
There exists $p \in (M \cap L) \setminus (N \cup N^\downarrow)$ such that $i_p \neq k_p$.
Here, $\mathsf{n}_p$ represents the number of $z_p$ satisfying $(x_p, z_p) \in R_{p,i_p}$ and $z_p = y_p$.
Hence, $\mathsf{n}_p = 0$ because $(x_p, y_p) \in R_{p,k_p}$ according to Lemma \ref{lem:i}.
\end{enumerate}
In all scenarios, we identified a value of $p \in P$ for which $\mathsf{n}_p = 0$.
\end{proof}

For convenience, we denote the empty sequence in the set $\mathbb{I}$ by $0$.
Let $R_i$ represent the set of all $i$-adjacent pairs in $X \times X$.

\begin{thm}[{cf.\ \cite[Theorem 3]{Bai}}]\label{thm:X}
With above notation, the pair $\mathfrak{X} = (X, \{R_i\}_{i \in \mathbb{I}})$ constitutes an association scheme.
Furthermore, if each $\mathfrak{X}_p$ is commutative, then $\mathfrak{X}$ is also commutative; similarly, if each $\mathfrak{X}_p$ is symmetric, then $\mathfrak{X}$ is likewise symmetric.
\end{thm}
\begin{proof}
It is straightforward to verify that conditions (R1) and (R2) are satisfied.
Lemma \ref{lem:i'} ensures that (R3) holds.
Lemmas \ref{lem:pijk1} and \ref{lem:pijk2} guarantee (R4).

Conditions (i) through (vi) in Lemma \ref{lem:pijk1} remain invariant under the interchange of $i$ and $j$.
Additionally, under the assumption of condition (iv), we observe that $\mathsf{k}_{p,i_p} = \mathsf{k}_{p,j_p}$ for $p \in (M \cap N) \setminus (L \cup L^\downarrow)$.
The number of $z$ as given by Lemma \ref{lem:pijk1} is unchanged when $i$ and $j$ are interchanged, provided that 
$\mathfrak{X}_p$ is commutative.
Consequently, if each $\mathfrak{X}_p$ is commutative, $\mathfrak{X}$ is also commutative.

If each $\mathfrak{X}_p$ is symmetric, Lemma \ref{lem:i'} implies (R6).
\end{proof}

We refer to the pair $\mathfrak{X} = (X, \{R_i\}_{i \in \mathbb{I}})$ as the generalized wreath product of the association schemes $\{\mathfrak{X}_p\}_{p \in P}$ over $P$.
It is noteworthy that \cite[Theorem 3]{Bai} assumed that each $\mathfrak{X}_p$ is symmetric, whereas we we do not assume commutativity in this context.
Nevertheless, since the proofs in \cite[Theorem 3]{Bai} do not depend on symmetry or commutativity, the results remain fundamentally unchanged.

\begin{prop}
With above notation, for $i = (i_p)_{p \in M} \in \mathbb{I}$, the $i$-th valency of $\mathfrak{X}$ is given by
\begin{equation}\label{k}
\mathsf{k}_i = \prod_{p \in M^\downarrow} |X_p| \times \prod_{p \in M} \mathsf{k}_{p,i_p}.
\end{equation}
\end{prop}
\begin{proof}
In Lemma \ref{lem:pijk1}, set $j = i'$ and $k = 0$.
\end{proof}

\section{Adjacency algebra}

Fix a total order on $P$.
Thereafter, unless otherwise specified, any Kronecker product $\bigotimes_{p \in P} W_p$ over the entire set $P$ will be taken according to this total order.
For each $p \in P$, let $\{A_{p,i}\}_{i=0}^{d_p}$ denote the adjacency matrices of $\mathfrak{X}_p$.
We abbreviate $I_p = I_{X_p}$ and $J_p = J_{X_p}$.
By Lemma \ref{lem:i}, the adjacency matrices of $\mathfrak{X}$ are given by
\begin{align}\label{A}
A_i = \bigotimes_{p \in M^\downarrow} J_p \otimes \bigotimes_{p \in M} A_{p,i_p} \otimes \bigotimes_{p \not\in M \cup M^\downarrow} I_p,
&&
i = (i_p)_{p \in M} \in \mathbb{I}.
\end{align}
Let $\mathcal{A}_p$ denote the adjacency algebra of $\mathfrak{X}_p$ for $p \in P$,
and let $\mathcal{A}$ denote the adjacency algebra of $\mathfrak{X}$.
Our objective in this section is to provide an expression for $\mathcal{A}$ in terms of $\mathcal{A}_p$ for $p \in P$.
In this paper, by a partition of $P$, we refer to a collection of mutually disjoint sets whose union is $P$,
where each element is not necessarily non-empty.

\begin{lem}\label{lem:Adj}
Let $D$, $M$ and $U$ denote a down-set, an anti-chain, and an up-set in $P$, respectively.
Assume $\{D,M,U\}$ forms a partition of $P$. 
For $i = (i_p)_{p \in M} \in \mathbb{I}$, we have
\[
\bigotimes_{p \in D} J_p \otimes \bigotimes_{p \in M} A_{p,i_p} \otimes \bigotimes_{p \in U} I_p \in \mathcal{A}.
\]
\end{lem}
\begin{proof}
Since $U$ is an up-set, we have $M^\downarrow \subset D$.
Let $Q = D \setminus M^\downarrow$.
If $Q = \emptyset$, then the matrix in the assertion coincides with the adjacency matrix \eqref{A}.
We may assume $Q \neq \emptyset$.
We define $\mathbb{I}(Q) = \{ (j_p)_{p \in N} \in \mathbb{I} \mid N \subset Q\}$.
For each $j = (j_p)_{p \in N} \in \mathbb{I}(Q)$,
we define $A_j(Q)$ analogously to $A_j$ as in \eqref{A} but with indices restricted to $Q$.
More precisely, 
\[
A_j(Q) = \bigotimes_{p \in Q \cap N^\downarrow} J_p \otimes \bigotimes_{p \in N} A_{p,j_p} \otimes \bigotimes_{p \in Q \setminus (N^\downarrow \cup N)} I_p.
\]
Then $\{A_j(Q)\}_{j \in \mathbb{I}(Q)}$ constitutes the complete set of adjacency matrices of the generalized wreath product of $\{\mathfrak{X}_p\}_{p \in Q}$ over $Q$.
In particular, the sum of all $A_j(Q)$ is the all-ones matrix.
For $j \in \mathbb{I}(Q)$, observe that $i \cup j \in \mathbb{I}$ and
\[
\sum_{j \in \mathbb{I}(Q)} A_{i \cup j} = \bigotimes_{p \in M^\downarrow} J_p \otimes \bigotimes_{p \in M} A_{p,i_p} \otimes \bigotimes_{p \in U} I_p \otimes \left(\sum_{j \in \mathbb{I}(Q)} A_j(Q)\right) = \bigotimes_{p \in D} J_p \otimes \bigotimes_{p \in M} A_{p,i_p} \otimes \bigotimes_{p \in U} I_p.
\]
The result follows.
\end{proof}

\begin{nota}\label{Delta}
Let $\Delta$ denote the set of ordered triples $(D,M,U)$ consisting of a down-set $D$, a non-empty anti-chain $M$ and an up-set $U$ such that $\{D,M,U\}$ forms a partition of $P$.
\end{nota}

\begin{lem}\label{lem:B}
Referring to Notation \ref{Delta},
for $(D,M,U) \in \Delta$, we have
\[
\bigotimes_{p \in D} J_p \otimes \bigotimes_{p \in M} \mathcal{A}_p \otimes \bigotimes_{p \in U} I_p \subset \mathcal{A}.
\]
\end{lem}
\begin{proof}
Take $0 \le i_p \le d_p$ for $p \in M$.
Since $\{A_{p,i}\}_{i=0}^{d_p}$ forms a basis for $\mathcal{A}_p$,
it suffices to show 
\[
A = \bigotimes_{p \in D} J_p \otimes \bigotimes_{p \in M} A_{p,i_p} \otimes \bigotimes_{p \in U} I_p \subset \mathcal{A}.
\]
Define $M_0 = \{p \in M \mid i_p = 0\}$, $M' = M \setminus M_0$ and $U' = U \cup M_0$.
Since $A_{p,0} = I_p$, the matrix $A$ can be expressed as
\[
A = \bigotimes_{p \in D} J_p \otimes \bigotimes_{p \in M'} A_{p,i_p} \otimes \bigotimes_{p \in U'} I_p.
\]
Since $M'$ is an anti-chain, $U'$ is an up-set, $\{D, M', U'\}$ is a partition of $P$ and $(i_p)_{p \in M'} \in \mathbb{I}$, $A$ is in the the adjacency algebra of $\mathfrak{X}$ by Lemma \ref{lem:Adj}.
\end{proof}
\begin{prop}\label{prop:Adj}
The adjacency algebra of $\mathfrak{X}$ is expressed as
\[
\mathcal{A} = \sum_{(D,M,U) \in \Delta} \left( \bigotimes_{p \in D} J_p \otimes \bigotimes_{p \in M} \mathcal{A}_p \otimes \bigotimes_{p \in U} I_p \right).
\]
\end{prop}
\begin{proof}
By Lemma \ref{lem:B}, each summand in the right-hand side of the desired equality is a subspace of $\mathcal{A}$.
In particular, the sum is also a subspace of $\mathcal{A}$.
It is remain to demonstrate the reverse inclusion.
Take $i = (i_p)_{p \in M} \in \mathbb{I}$ and 
let $D = M^\downarrow$ and $U = P \setminus (M \cup M^\downarrow)$.
Since $A_{p,i_p} \in \mathcal{A}_p$ for $p \in M$,
the adjacency matrix $A_i$ defined in \eqref{A} belongs to the summand corresponding to $(D,M,U)$.
Since $\{A_i\}_{i \in \mathbb{I}}$ forms a basis for $\mathcal{A}$,
the reverse inclusion follows.
\end{proof}

\begin{cor}\label{cor:Adj}
Let $Q$ be a non-empty subposet of $P$.
Let $\mathfrak{X}_Q$ denote the generalized wreath product of the association schemes $\{\mathfrak{X}_p\}_{p \in Q}$ over $Q$.
Let $B_p \in \mathcal{A}_p$ for $p \in P$.
If $\bigotimes_{p \in P} B_p \in \mathcal{A}$, then
$\bigotimes_{p \in Q} B_p$ belongs to the adjacency algebra of $\mathfrak{X}_Q$.
\end{cor}
\begin{proof}
Let $B = \bigotimes_{p \in P} B_p$ and $B_Q = \bigotimes_{p \in Q} B_p$.
By Proposition \ref{prop:Adj},
since $\{A_{p,i}\}_{i=0}^{d_p}$ constitutes a basis for $\mathcal{A}_p$,
we may assume 
\[
B = \bigotimes_{p \in D} J_p \otimes \bigotimes_{p \in M} A_{p,i_p} \otimes \bigotimes_{p \in U} I_p
\]
for some $(D,M,U) \in \Delta$ and $0 \le i_p \le d_p$.
Thus,
\[
B_Q = \bigotimes_{p \in D \cap Q} J_p \otimes \bigotimes_{p \in M \cap Q} A_{p,i_p} \otimes \bigotimes_{p \in U \cap Q} I_p.
\]
Since $D$ is a down-set in $P$, $D \cap Q$ is also a down-set in $Q$.
Similarly, $M \cap Q$ is an anti-chain in $Q$ and $U \cap Q$ is an up-set in $Q$.
Furthermore, $\{D \cap Q, M \cap Q, U \cap Q\}$ forms a partition of $Q$.
If $M \cap Q = \emptyset$, $B_Q$ belongs to the adjacency algebra of $\mathfrak{X}_Q$ by Lemma \ref{lem:Adj}.
If $M \cap Q \neq \emptyset$, $B_Q$ is in the adjacency algebra of $\mathfrak{X}_Q$ by Lemma \ref{lem:B}.
\end{proof}

\section{Irreducible representations of adjacency algebra}
For each $p \in P$, since the adjacency algebra $\mathcal{A}_p$ of $\mathfrak{X}_p$ is semi-simple, it can be decomposed as a direct sum of central simple algebras: 
\begin{equation}\label{Ap}
\mathcal{A}_p = \bigoplus_{j=0}^{l_p} \mathcal{A}_{p,j}.
\end{equation}
As customary, we assume $\mathcal{A}_{p,0} = \langle J_p \rangle$.
In this section, we decompose the adjacency algebra $\mathcal{A}$ of $\mathfrak{X}$ into central simple algebras.
According to the decomposition \eqref{Ap}, we define
\[
\mathbb{J} = \{ (j_p)_{p \in M} \mid \text{$M$ is an anti-chain in $P$}, 1 \le j_p \le l_p \; (p \in M)\}.
\]
We denote the empty sequence in the set $\mathbb{J}$ by $0$, as we did for $\mathbb{I}$.
Then we set
\begin{align}\label{Aj}
\mathcal{A}_j = 
\bigotimes_{p \not\in M \cup M^\uparrow} J_p \otimes \bigotimes_{p \in M} \mathcal{A}_{p,j_p} \otimes \bigotimes_{p \in M^\uparrow} I_p,
&&
j = (j_p)_{p \in M} \in \mathbb{J} \setminus \{0\}.
\end{align}
For convenience, we set $\mathcal{A}_0 = \langle J_X \rangle$.

\begin{lem}\label{lem:AjinA}
With above notation, for $j \in \mathbb{J}$, we have $\mathcal{A}_j \subset \mathcal{A}$.
\end{lem}
\begin{proof}
If $j = 0$, since $J_X \in \mathcal{A}$, the assertion holds.
Assume $j \neq 0$ and let $j = (j_p)_{p \in M}$.
Define $D = P \setminus (M \cup M^\uparrow)$ and $U = M^\uparrow$.
Then $(D,M,U) \in \Delta$ as per Notation \ref{Delta}.
By Lemma \ref{lem:B},
since $\mathcal{A}_{p,j_p} \subset \mathcal{A}_p$, it follows that $\mathcal{A}_j \subset \mathcal{A}$.
\end{proof}

\begin{lem}\label{lem:ideal}
With above notation, for $j \in \mathbb{J}$, the subspace $\mathcal{A}_j$ constitutes a two-sided ideal of $\mathcal{A}$.
\end{lem}
\begin{proof}
If $j = 0$, the assertion is evident since $J_X \in \mathcal{A}$.
Assume $j \neq 0$.
We first demonstrate $\mathcal{A}_j$ is a left ideal of $\mathcal{A}$.
Since $\{A_i\}_{i \in \mathbb{I}}$ forms a basis for $\mathcal{A}$,
it suffices to show $A_i\mathcal{A}_j \subset \mathcal{A}_j$ for $i \in \mathbb{I}$.
Let $j = (j_p)_{p \in M}$ and take $i = (i_p)_{p \in N} \in \mathbb{I}$.

If there exists $p \in M \cap N^\downarrow$, since $j_p \neq 0$, the $p$-th factor of $A_i\mathcal{A}_j$ is $J_p\mathcal{A}_{p,j_p} = 0$.
If there exists $p \in M^\uparrow \cap (N \cup N^\downarrow)$, 
there exists $q \in M \cap N^\downarrow$ such that $q < p$.
Therefore, we have $A_i\mathcal{A}_j = 0 \subset \mathcal{A}_j$ if $M \cap N^\downarrow \neq \emptyset$ or $M^\uparrow \cap (N \cup N^\downarrow) \neq \emptyset$.

Consider the case where $M \cap N^\downarrow = M^\uparrow \cap (N \cup N^\downarrow) = \emptyset$.
If $p \in M^\uparrow$, then $p \not\in N \cup N^\downarrow$.
This implies that if the $p$-th factor of $\mathcal{A}_j$ is $I_p$, then the $p$-th factor of $A_i$ is also $I_p$.
The remianing factor of $\mathcal{A}_j$ is either $J_p$ or $\mathcal{A}_{p,j_p}$, both of which are two-sided ideals of $\mathcal{A}_p$.
Thus, $A_i\mathcal{A}_j \subset \mathcal{A}_j$ in this scenario.
This establishes that $\mathcal{A}_j$ is a left ideal of $\mathcal{A}$.

As similar argument demonstrates that $\mathcal{A}_j$ is also a right ideal of $\mathcal{A}$.
Hence, the result follows.
\end{proof}

\begin{lem}\label{lem:orthogonal}
With above notation, the subspaces $\{\mathcal{A}_j\}_{j \in \mathbb{J}}$ are mutually orthogonal.
\end{lem}
\begin{proof}
Let $i = (i_p)_{p \in M}, j = (j_p)_{p \in N} \in \mathbb{J}$ with $i \neq j$ and demonstrate that $\mathcal{A}_i\mathcal{A}_j = 0$.
It suffices to show that there exists $p \in P$ such that the $p$-th factor of $J_X\mathcal{A}_j$ is zero.

If $i = 0$ and $j \neq 0$, since there exists $p \in N$ with $j_p \neq 0$, the $p$-th factor of $J_X\mathcal{A}_j$ is $J_p\mathcal{A}_{p,j_p} = 0$.
Similarly, the assertion holds if $i \neq 0$ and $j = 0$.

Assume $i \neq 0$, $j \neq 0$ and $M \not\subset N$.
In this scenario, there exists $p \in M$ such that $p \not\in N$.
Then either $p \not\in N \cup N^\uparrow$ or $p \in N^\uparrow$.
If $p \not\in N \cup N^\uparrow$ the $p$-th factor of $\mathcal{A}_i\mathcal{A}_j$ is $\mathcal{A}_{p,i_p}J_p = 0$ since $i_p \neq 0$.
If $p \in N^\uparrow$, there exists $q \in M \setminus (N \cup N^\uparrow)$ such that $q < p$, which reduces to a previous case.
Therefore, the assertion holds if $i \neq 0$, $j \neq 0$ and $M \not\subset N$.
Similarly, the assertion holds if $i \neq 0$, $j \neq 0$ and $N \not\subset M$.

Consider the case $M = N$.
Since $i \neq j$, there exists $p \in M = N$ such that $i_p \neq j_p$.
Thus the $p$-th factor of $\mathcal{A}_i\mathcal{A}_j$ is $\mathcal{A}_{p,i_p}\mathcal{A}_{p,j_p} = 0$.
\end{proof}

\begin{thm}\label{thm:A}
The adjacency algebra $\mathcal{A}$ of $\mathfrak{X}$ is represented as a direct sum of central simple algebras:
\[
\mathcal{A} = \bigoplus_{j \in \mathbb{J}} \mathcal{A}_j,
\]
where $\mathcal{A}_j$ is specified in \eqref{Aj}.
\end{thm}
\begin{proof}
First, we demonstrate that each $\mathcal{A}_j$ is a central simple algebra.
If $j=0$, since $\mathcal{A}_0$ has dimension $1$, it is a central simple algebra.
For $j = (j_p)_{p \in M} \in \mathbb{J}$ with $j \neq 0$,
by definition \eqref{Aj}, we have as algebras
\[
\mathcal{A}_j \cong \bigotimes_{p \in M} \mathcal{A}_{p,j_p}
\]
as algebras. 
Hence, $\mathcal{A}_j$ is a central simple algebra as it is a tensor product of central simple algebras.
Next, we verify the desired equality.
By Lemma \ref{lem:AjinA}, it follows that
\[
\sum_{j \in \mathbb{J}} \mathcal{A}_j \subset \mathcal{A}.
\]
According to Lemmas \ref{lem:ideal} and \ref{lem:orthogonal}, the sum is direct.
To establish the reverse inclusion, we compare the dimensions of both sides:
\[
\dim \left(\sum_{j \in \mathbb{J}} \mathcal{A}_j \right) 
= \sum_{j \in \mathbb{J}} \dim \left(\mathcal{A}_j \right) 
= \sum_{M} \prod_{p \in M} \sum_{j_p = 1}^{l_p} \left( \dim \mathcal{A}_{p,j_p} \right) = \sum_{M}  \prod_{p \in M} d_p = |\mathbb{I}| = \dim \mathcal{A},
\]
where the third and fifth sums iterate over all anti-chains $M$ in $P$.
\end{proof}

\begin{thm}\label{thm:E}
For each $p \in P$, let $\{E_{p,j}\}_{j=0}^{l_p}$ denote the primitive central idempotents of the adjacency algebra of $\mathfrak{X}_p$
with $E_{p,0} = |X_p|^{-1}J_p$.
Then, the primitive central idempotents of the adjacency algebra of $\mathfrak{X}$ are given by
\begin{align*}
\bigotimes_{p \not\in M \cup M^\uparrow} E_{p,0} \otimes \bigotimes_{p \in M} E_{p,j_p} \otimes \bigotimes_{p \in M^\uparrow} I_p,
&&
j = (j_p)_{p \in M} \in \mathbb{J}.
\end{align*}
\end{thm}
\begin{proof}
By Theorem \ref{thm:A}, every primitive central idempotent is the identity in $\mathcal{A}_j$ for some $j \in \mathbb{J}$.
The result follows from definition \eqref{Aj} and the definition of $\{E_{p,j}\}_{j=0}^{l_p}$.
\end{proof}

\section{Dual adjacency algebra}
From now on until end of this paper, we fix an element $x = (x_p)_{p \in P} \in X$.
For each $p \in P$, let $\{E_{p,i_p}^*(x_p)\}_{i_p=0}^{d_p}$ denote the dual primitive idempotents  of $\mathfrak{X}_p$ with respect to $x_p$.
We abbreviate $E_{p,i_p}^* = E_{p,i_p}^*(x_p)$ for $p \in P$.
According to \eqref{A}, the dual primitive idempotents of $\mathfrak{X}$ with respect to $x$ are given by
\begin{align}\label{E*}
E_i^*(x) = \bigotimes_{p \in M^\downarrow} I_p  \otimes \bigotimes_{p \in M} E_{p,i_p}^* \otimes \bigotimes_{p \not\in M \cup M^\downarrow} E_{p,0}^*,
&&
i = (i_p)_{p \in M} \in \mathbb{I}.
\end{align}
Let $\mathcal{A}_p^* = \mathcal{A}_p^*(x_p)$ denote the dual adjacency algebra of $\mathfrak{X}_p$ with respect to $x_p$ for $p \in P$
and
let $\mathcal{A}^*(x)$ denote the dual adjacency algebra of $\mathfrak{X}$ with respect to $x$.
Our objective in this section is to provide a precise representation of $\mathcal{A}^*(x)$ in terms of $\mathcal{A}_p^*$ for $p \in P$.

\begin{lem}\label{lem:B*}
Referring to Notation \ref{Delta},
for $(D,M,U) \in \Delta$, we have
\[
\bigotimes_{p \in D} I_p \otimes \bigotimes_{p \in M} \mathcal{A}_p^* \otimes \bigotimes_{p \in U} E_{p,0}^* \subset \mathcal{A}^*(x).
\]
\end{lem}
\begin{proof}
Analogous to Lemma \ref{lem:B}.
\end{proof}

\begin{prop}\label{prop:dAdj}
The dual adjacency algebra of $\mathfrak{X}$ with respect to $x$ is expressed as
\[
\mathcal{A}^*(x) = \sum_{(D,M,U) \in \Delta} \left( \bigotimes_{p \in D} I_p \otimes \bigotimes_{p \in M} \mathcal{A}_p^* \otimes \bigotimes_{p \in U} E_{p,0}^* \right).
\]
\end{prop}
\begin{proof}
Analogous to Proposition \ref{prop:Adj}.
\end{proof}

\begin{cor}\label{cor:dAdj}
Let $Q$ be a non-empty subposet of $P$.
Let $\mathfrak{X}_Q$ denote the generalized wreath product of the association schemes $\{\mathfrak{X}_p\}_{p \in Q}$ over $Q$.
Let $B_p^* \in \mathcal{A}_p^*$ for $p \in P$.
If $\bigotimes_{p \in P} B_p^* \in \mathcal{A}^*(x)$, then
$\bigotimes_{p \in Q} B_p^*$ belongs to the dual adjacency algebra of $\mathfrak{X}_Q$ with respect to $(x_p)_{p \in Q}$.
\end{cor}
\begin{proof}
Analogous to Corollary \ref{cor:Adj}.
\end{proof}

\section{Terwilliger algebra}

Let $\mathcal{T}_p = \mathcal{T}_p(x_p)$ denote the Terwilliger algebra of $\mathfrak{X}_p$ with respect to $x_p$ for $p \in P$
and
let $\mathcal{T}(x)$ denote the Terwilliger algebra of $\mathfrak{X}$ with respect to $x$.
Our objective in this section is to provide an explicit description of the generators of $\mathcal{T}(x)$ in terms of $\mathcal{A}_p$, $\mathcal{A}_p^*$, $\mathcal{T}_p$ for $p \in P$.

\begin{nota}\label{Omega}
Let $\Omega$ denote the set of ordered sextuples $(D,M,N,L,C,U)$ 
such that 
$(D, M \cup N, L \cup C \cup U), (D \cup M \cup C, N \cup L, U) \in \Delta$,
$C = M^\uparrow \cap L^\downarrow$, $M \subset L^\downarrow$ and $L \subset M^\uparrow$.
\end{nota}

\begin{lem}\label{lem:DQU}
Let $D$ be a down-set, $Q$ be a non-empty convex subset, $U$ be an up-set.
Assume $\{D,Q,U\}$ forms a partition of $P$.
Define
$M = \min(Q) \setminus \max(Q)$,
$N = \max(Q) \cap \min(Q)$,
$L = \max(Q) \setminus \min(Q)$,
$C = Q \setminus (\max(Q) \cup \min(Q))$.
Then $(D,M,N,L,C,U) \in \Omega$.
\end{lem}
\begin{proof}
Since it is the set of minimal (resp.\ maximal) elements, 
$M \cup N = \min(Q)$ (resp.\ $N \cup L = \max(Q)$) is an anti-chain.
Moreover, since $Q$ is non-empty, both $M \cup N$ and $N \cup L$ are non-empty.
Observe that $L \cup C = Q \setminus \min(Q)$.
Since $Q \cup U$ is an up-set,
if $p \in L \cup C$ and $p < q$, we have $q \in Q \cup U$ and $q \not\in \min(Q)$.
Since $U$ is an up-set, this implies that $L \cup C \cup U$ is an up-set.
Similarly, $D \cup M \cup C$ is a down-set.
By construction, we have $C = M^\uparrow \cap L^\downarrow$, $M \subset L^\downarrow$ and $L \subset M^\uparrow$.
The result follows.
\end{proof}

\begin{lem}\label{lem:S}
Referring to Notation \ref{Omega},
for $(D,M,N,L,C,U) \in \Omega$, we have
\[
\bigotimes_{p \in D} J_p
\otimes \bigotimes_{p \in M} \mathcal{A}_p
\otimes \bigotimes_{p \in N} \mathcal{T}_p
\otimes \bigotimes_{p \in L} \mathcal{A}_p^*
\otimes \bigotimes_{p \in C} I_p
\otimes \bigotimes_{p \in U} E_{p,0}^* \subset \mathcal{T}(x).
\]
\end{lem}
\begin{proof}
Let 
$D_1 = D$,
$M_1 = M \cup N$,
$U_1 = L \cup C \cup U$,
$D_2 = D \cup M \cup C$,
$M_2 = N \cup L$,
$U_2 = U$.
By the definition of $\Omega$,
we have $(D_1,M_1,U_1), (D_2,M_2,U_2) \in \Delta$.
By Lemmas \ref{lem:B} and \ref{lem:B*},
it follows that
\begin{align*}
\bigotimes_{p \in D_1} J_p \otimes \bigotimes_{p \in M_1} \mathcal{A}_p \otimes \bigotimes_{p \in U_1} I_p \subset \mathcal{A},
&&
\bigotimes_{p \in D_2} I_p \otimes \bigotimes_{p \in M_2} \mathcal{A}_p^* \otimes \bigotimes_{p \in U_2} E_{p,0}^* \subset \mathcal{A}^*(x).
\end{align*}
Consider the subspace $\mathcal{S}(x)$ of $\mathcal{T}(x)$ generated by these two subspaces.
Since $\mathcal{T}_p$ is generated by $\mathcal{A}_p$ and $\mathcal{A}_p^*$, we have
\[
\mathcal{S}(x) = \bigotimes_{p \in D} J_p
\otimes \bigotimes_{p \in M} \mathcal{A}_p
\otimes \bigotimes_{p \in N} \mathcal{T}_p
\otimes \bigotimes_{p \in L} \mathcal{A}_p^*
\otimes \bigotimes_{p \in C} I_p
\otimes \bigotimes_{p \in U} E_{p,0}^*.
\]
Since $\mathcal{S}(x)$ is the subspace of $\mathcal{T}(x)$,
the result follows.
\end{proof}

\begin{prop}\label{prop:T}
The Terwilliger algebra of $\mathfrak{X}$ with respect to $x$ is generated by
\begin{align*}
\bigotimes_{p \in D} J_p
\otimes \bigotimes_{p \in M} \mathcal{A}_p
\otimes \bigotimes_{p \in N} \mathcal{T}_p
\otimes \bigotimes_{p \in L} \mathcal{A}_p^*
\otimes \bigotimes_{p \in C} I_p
\otimes \bigotimes_{p \in U} E_{p,0}^*,
&&
(D,M,N,L,C,U) \in \Omega.
\end{align*}
\end{prop}
\begin{proof}
Since $\mathcal{T}(x)$ is generated by $\mathcal{A}$ and $\mathcal{A}^*(x)$,
it suffices to demonstrate the existence of $(D_1,M_1,N_1,L_1,C_1,U_1)$, $(D_2,M_2,N_2,L_2,C_2,U_2) \in \Omega$ such that 
\begin{equation}
\label{AinS}
\mathcal{A} \subset \bigotimes_{p \in D_1} J_p
\otimes \bigotimes_{p \in M_1} \mathcal{A}_p
\otimes \bigotimes_{p \in N_1} \mathcal{T}_p
\otimes \bigotimes_{p \in L_1} \mathcal{A}_p^*
\otimes \bigotimes_{p \in C_1} I_p
\otimes \bigotimes_{p \in U_1} E_{p,0}^*,
\end{equation}
and
\begin{equation}
\label{A*inS}
\mathcal{A}^*(x) \subset \bigotimes_{p \in D_2} J_p
\otimes \bigotimes_{p \in M_2} \mathcal{A}_p
\otimes \bigotimes_{p \in N_2} \mathcal{T}_p
\otimes \bigotimes_{p \in L_2} \mathcal{A}_p^*
\otimes \bigotimes_{p \in C_2} I_p
\otimes \bigotimes_{p \in U_2} E_{p,0}^*.
\end{equation}
Let $(D,M,U) \in \Delta$ and define $Q = M \cup U$.
Then set
$D_1 = D$,
$M_1 = \min(Q) \setminus \max(Q)$,
$N_1 = \max(Q) \cap \min(Q)$,
$L_1 = \max(Q) \setminus \min(Q)$,
$C_1 = Q \setminus (\max(Q) \cup \min(Q))$,
$U_1 = \emptyset$.
By Lemma \ref{lem:DQU}, we have $(D_1,M_1,N_1,L_1,C_1,U_1) \in \Omega$.
If $p \in M \subset Q$ and $q < p$, then $q \in D$, implying $q \not\in Q$.
Consequently, $p \in \min(Q) = M_1 \cup N_1$, thus $M \subset M_1 \cup N_1$.
Since $\mathcal{A}_p \subset \mathcal{T}_p$ and $I_p \in \mathcal{A}_p \cap \mathcal{A}_p^* \subset \mathcal{T}_p$, we obtain
\[
\bigotimes_{p \in D} J_p \otimes \bigotimes_{p \in M} \mathcal{A}_p \otimes \bigotimes_{p \in U} I_p \subset \bigotimes_{p \in D_1} J_p
\otimes \bigotimes_{p \in M_1} \mathcal{A}_p
\otimes \bigotimes_{p \in N_1} \mathcal{T}_p
\otimes \bigotimes_{p \in L_1} \mathcal{A}_p^*
\otimes \bigotimes_{p \in C_1} I_p
\otimes \bigotimes_{p \in U_1} E_{p,0}^*.
\]
Therefore, by Proposition \ref{prop:Adj}, inclusion \eqref{AinS} follows.
Similarly, by Proposition \ref{prop:dAdj}, inclusion \eqref{A*inS} is established.
Hence, the result follows.
\end{proof}

\begin{cor}\label{cor:T}
Let $Q$ be a non-empty subposet of $P$.
Let $\mathfrak{X}_Q$ denote the generalized wreath product of the association schemes $\{\mathfrak{X}_p\}_{p \in Q}$ over $Q$.
Let $S_p \in \mathcal{T}_p$ for $p \in P$.
If $\bigotimes_{p \in P} S_p \in \mathcal{T}(x)$, then
$\bigotimes_{p \in Q} S_p$ belongs to the Terwilliger algebra of $\mathfrak{X}_Q$ with respect to $(x_p)_{p \in Q}$.
\end{cor}
\begin{proof}
This follows directly from Corollaries \ref{cor:Adj} and \ref{cor:dAdj}.
\end{proof}

\section{Principal two-sided ideal generated by the all-ones matrix}

Let $\mathcal{J}_p = \mathcal{J}_p(x_p)$ represent the principal two-sided ideal of $\mathcal{T}_p$ generated by the all-ones matrix $J_p$ for $p \in P$,
and
let $\mathcal{J}(x)$ denote the principal two-sided ideal of $\mathcal{T}(x)$ generated by the all-ones matrix $J_X$.
The objective of this section is to determine the structure of $\mathcal{J}(x)$.
By Lemma \ref{lem}, it follows that $\mathcal{J}_p = \mathcal{A}^*_pJ_p\mathcal{A}^*_p$ for $p \in P$, and similarly, $\mathcal{J}(x) = \mathcal{A}^*(x)J_X\mathcal{A}^*(x)$.

\begin{lem}\label{lem:K}
Referring to Notation \ref{Delta},
for $(D,M,U) \in \Delta$, we have
\[
\bigotimes_{p \in D} J_p \otimes \bigotimes_{p \in M} \mathcal{J}_p \otimes \bigotimes_{p \in U} E_{p,0}^* \subset \mathcal{J}(x).
\]
\end{lem}
\begin{proof}
By Lemma \ref{lem:B*}, we obtain
\[
\bigotimes_{p \in D} I_p \otimes \bigotimes_{p \in M} \mathcal{A}_p^* \otimes \bigotimes_{p \in U} E_{p,0}^* \subset \mathcal{A}^*(x).
\]
Given that $I_pJ_pI_p = J_p$, $\mathcal{A}_p^*J_p\mathcal{A}_p^* = \mathcal{J}_p$ and $E_{p,0}^*J_pE_{p,0}^* = E_{p,0}^*$,
the matrix below is derived by multiplying the aforementioned subspace by the all-ones matrix $J_X$ on both sides:
\[
\bigotimes_{p \in D} J_p \otimes \bigotimes_{p \in M} \mathcal{J}_p \otimes \bigotimes_{p \in U} E_{p,0}^*.
\]
By construction, this subspace lies within $\mathcal{J}(x)$.
\end{proof}
\begin{thm}\label{thm:J}
The principal two-sided ideal of $\mathcal{T}(x)$ generated by the all-ones matrix $J_X$ is expressed as $\mathcal{J}(x) = \mathcal{K}(x)\mathcal{K}(x)$, where
\[
\mathcal{K}(x) = \sum_{(D,M,U) \in \Delta} \left( \bigotimes_{p \in D} J_p \otimes \bigotimes_{p \in M} \mathcal{J}_p \otimes \bigotimes_{p \in U} E_{p,0}^* \right).
\]
\end{thm}
\begin{proof}
By Lemma \ref{lem:K}, we have $\mathcal{K}(x)\mathcal{K}(x) \subset \mathcal{J}(x)$.
It remains to demonstrate the reverse inclusion.
Referring to Lemma \ref{lem:B*},
for $\delta = (D,M,U) \in \Delta$, we set
\[
\mathcal{B}^*_{\delta}(x) = \bigotimes_{p \in D} I_p \otimes \bigotimes_{p \in M} \mathcal{A}_p^* \otimes \bigotimes_{p \in U} E_{p,0}^* \subset \mathcal{A}^*(x).
\]
By Lemma \ref{lem:K}, it is noted that $\mathcal{B}^*_{\delta}(x) J_X \mathcal{B}^*_{\delta}(x) \subset \mathcal{K}(x)$ for $\delta \in \Delta$.
Since $\mathcal{J}(x) = \mathcal{A}^*(x)J_X\mathcal{A}^*(x)$,
by Proposition \ref{prop:dAdj},
it suffices to show that $\mathcal{B}^*_{\delta_1}(x) J_X \mathcal{B}^*_{\delta_2}(x) \subset \mathcal{K}(x)\mathcal{K}(x)$ for $\delta_1, \delta_2 \in \Delta$.
Given that $J_X \mathcal{T}(x) J_X = \langle J_X \rangle$,
we obtain $\mathcal{B}^*_{\delta_1}(x) J_X \mathcal{B}^*_{\delta_2}(x) = \mathcal{B}^*_{\delta_1}(x) J_X \mathcal{B}^*_{\delta_1}(x)\mathcal{B}^*_{\delta_2}(x) J_X \mathcal{B}^*_{\delta_2}(x) \subset \mathcal{K}(x)\mathcal{K}(x)$.
Thu, the result follows.
\end{proof}

\begin{thm}\label{thm:eps}
For each $p \in P$, let $\epsilon_p = \epsilon_p(x_p)$ denote the principal idempotent of $\mathcal{T}_p$.
Then the principal idempotent of $\mathcal{T}(x)$ is expressed as
\[
\sum_{M}
\left(
\bigotimes_{p \in M^\downarrow} E_{p,0}
\otimes \bigotimes_{p \in M} (\epsilon_p - E_{p,0}^*)
\otimes \bigotimes_{p \not\in M \cup M^\downarrow} E_{p,0}^*
\right),
\]
where the summation is taken over all anti-chains $M$ in $P$.
\end{thm}
\begin{proof}
The principal idempotent is expressed as
$\sum_{i \in \mathbb{I}} \mathsf{k}_i^{-1} E_i^* J_X E_i^*$.
By utilizing the identities \eqref{k} and \eqref{E*}, we obtain
\begin{align*}
\sum_{i \in \mathbb{I}} \frac{1}{\mathsf{k}_i} E_i^* J_X E_i^*
&= \sum_{(i_p)_{p \in M} \in \mathbb{I}} \left(\bigotimes_{p \in M^\downarrow} \frac{1}{|X_p|}J_p  \otimes \bigotimes_{p \in M} \frac{1}{\mathsf{k}_{p,i}}E_{p,i_p}^*J_pE_{p,i_p}^* \otimes \bigotimes_{p \not\in M \cup M^\downarrow} E_{p,0}^*,\right) \\
&= \sum_{M} \left(\bigotimes_{p \in M^\downarrow} \frac{1}{|X_p|}J_p  \otimes \bigotimes_{p \in M} \sum_{i_p = 1}^{d_p}\frac{1}{\mathsf{k}_{p,i}}E_{p,i_p}^*J_pE_{p,i_p}^* \otimes \bigotimes_{p \not\in M \cup M^\downarrow} E_{p,0}^*,\right),
\end{align*}
where the final summation is taken over all anti-chains $M$ in $P$.
Since $E_{p,0} = |X_p|^{-1}J_p$, $\epsilon_p = \sum_{i_p = 0}^{d_p}\mathsf{k}_{p,i}^{-1}E_{p,i_p}^*J_pE_{p,i_p}^*$,
$\mathsf{k}_{p,0}^{-1}E_{p,0}^*J_pE_{p,0}^* = E_{p,0}^*$,
the desired result follows.
\end{proof}

\section{Irreducible representations of terwilliger algebra}

For each $p \in P$,
since the Terwilliger algebra $\mathcal{T}_p$ of $\mathfrak{X}_p$ is semi-simple, it can be decomposed into a direct sum of central simple algebras: 
\begin{equation}\label{Tp}
\mathcal{T}_p = \bigoplus_{k=0}^{t_p} \mathcal{T}_{p,k}(x_p).
\end{equation}
By convention, we set $\mathcal{T}_{p,0}(x_p) = \mathcal{J}_p$.
We denote $\mathcal{T}_{p,k} = \mathcal{T}_{p,k}(x_p)$.
In this section, we decompose the Terwilliger algebra $\mathcal{T}(x)$ of $\mathfrak{X}$ into central simple algebras.
According to the decomposition \eqref{Tp}, we define
\[
\mathbb{K} = \{ (k_p)_{p \in M} \mid \text{$M$ is an anti-chain in $P$}, 1 \le k_p \le t_p \; (p \in M)\}.
\]
As with $\mathbb{I}$ and $\mathbb{J}$, we denote the empty sequence in the set $\mathbb{K}$ by $0$.

\begin{nota}\label{Lambda}
Let $\Lambda$ denote the set of ordered triple $(M,N,L)$ such that 
$M \subset L^\downarrow$, $L \subset M^\uparrow$
and that
$M \cup N$, $N \cup L$ are anti-chains.
\end{nota}

\begin{lem}\label{lem:bij} 
For $(M,N,L) \in \Lambda$, there exists a bijection between the following two sets:
\begin{enumerate}
\item The set of ordered triples $(D,H,U)$ such that $N \cap H = \emptyset$ and $(D,M,N \cup H,L,M^\uparrow \cap L^\downarrow,U) \in \Omega$.
\item The set of ordered triples $(D,H,U)$ comprising a down-set $D$ in $Q$, an anti-chain $H$ in $Q$, an up-set $U$ in $Q$ such that $\{D,H,U\}$ forms a partition of $Q$,
where $Q = P \setminus (M^\uparrow \cup N^\uparrow \cup N \cup N^\downarrow \cup L^\downarrow)$.
If $Q = \emptyset$, this set is interpreted as containing only $(\emptyset, \emptyset, \emptyset)$.
\end{enumerate}
In patricular, a bijection from set (i) to set (ii) is given by $\phi(D,H,U) = (D \cap Q,H,U \cap Q)$.
Moreover, both sets are non-empty for any $(M,N,L) \in \Lambda$.
\end{lem}
\begin{proof}
Fix $(M,N,L) \in \Lambda$ and let $Q = P \setminus (M^\uparrow \cup N^\uparrow \cup N \cup N^\downarrow \cup L^\downarrow)$.

Consider an element $(D_1,H_1,U_1)$ from set (i).
If $p \in D_1 \cap Q$, $q \in Q$ with $q < p$, we have $q \in D_1$ because $D_1$ is a down-set in $P$.
Hence, $D_1 \cap Q$ is a down-set in $Q$.
Similarly, $U_1 \cap Q$ is an up-set in $Q$.
If $p \in H_1$, since $N \cap H_1 = \emptyset$ and $(D_1,M,N \cup H_1,L,M^\uparrow \cap L^\downarrow,U_1) \in \Omega$, 
it follows that $p \in Q$.
Thus, $H_1 \cap Q = H_1$.
By the definition of $Q$, we have $M \cap Q = N \cap Q = L \cap Q = M^\uparrow \cap L^\downarrow \cap Q = \emptyset$.
This demonstrates that $\{D_1 \cap Q,H_1,U_1 \cap Q\}$ constitutes a partition of $Q$.
Therefore, the ordered triple $(D_1 \cap Q,H_1,U_1 \cap Q)$ belongs to set (ii).
Consequently, we define a map $\phi$ from set (i) to set (ii) by $\phi(D_1,H_1,U_1) = (D_1 \cap Q,H_1,U_1 \cap Q)$.

It remains to demonstrate that the map $\phi$ is a bijection.
Define $D' = N^\downarrow \cup (L^\downarrow \setminus (M \cup M^\uparrow))$
and 
$U' = (M^\uparrow \cup N^\uparrow) \setminus (L \cup L^\downarrow)$.
Conciser an ordered triple $(D_1,H_1,U_1)$ from set (i).
Since $(D_1,M,N \cup H_1,L,M^\uparrow \cap L^\downarrow,U_1) \in \Omega$,
it follows that $D' \subset D_1$.
If $p \in D_1 \setminus D'$, we have $p \not\in M \cup N$ since $D_1 \cap (M \cup N) = \emptyset$,
and $p \not\in N^\downarrow$ since $N^\downarrow \subset D'$.
Additionally, $p \not\in M^\uparrow \cup N^\uparrow$ since $D_1$ is a down-set in $P$
If $p \in (D_1 \setminus D') \cap L^\downarrow$, then by the definition of $D'$, $p \in M \cup M^\uparrow$,
which is a contradiction since $D_1 \cap M = \emptyset$ and $D_1$ is a down-set in $P$.
Thus, we have established that $D_1 \setminus D' \subset Q$.
Since $D' \cap Q = \emptyset$, it follows that $D_1 \setminus D' = D_1 \cap Q$.
Similarly, we have $U' \subset U_1$ and $U_1 \setminus U' = U_1 \cap Q$.
Consequently, 
$\phi(D_1,H_1,U_1) = (D_1 \setminus D',H_1,U_1 \setminus U')$.
Given that $D' \subset D_1$ and $U' \subset U_1$, the map $\phi$ is indeed a bijection.

Finally, since set (ii) includes both $(Q,\emptyset, \emptyset)$ and $(\emptyset, \emptyset, Q)$, it follows that both sets are non-empty.
\end{proof}

We define
\[
\mathbb{T} = \left\{\left((i_p)_{p \in L}, (j_p)_{p \in M}, (k_p)_{p \in N}\right) \in \mathbb{I} \times \mathbb{J} \times \mathbb{K} \mid (M,N,L) \in \Lambda\right\}.
\]
For $(i,j,k) \in \mathbb{T}$ with
$i = (i_p)_{p \in L}$,
$j = (j_p)_{p \in M}$, and
$k = (k_p)_{p \in N}$,
let $\mathcal{T}_{i,j,k}(x)$ denote the algebra generated by
\begin{align}\label{Tijk}
\bigotimes_{p \in D} J_p
\otimes \bigotimes_{p \in M} \mathcal{A}_{p,j_p}
\otimes \bigotimes_{p \in N} \mathcal{T}_{p,k_p}
\otimes \bigotimes_{p \in H} \mathcal{J}_p
\otimes \bigotimes_{p \in L} E_{p,i_p}^*
\otimes \bigotimes_{p \in M^\uparrow \cap L^\downarrow} I_p
\otimes \bigotimes_{p \in U} E_{p,0}^*,
\end{align}
for $(D,H,U)$ in the set Lemma \ref{lem:bij} (i) with respect to $(M,N,L)$.

\begin{lem}\label{lem:Q=empty}
With above notation, let $(i,j,k) \in \mathbb{T}$ with
$i = (i_p)_{p \in L}$,
$j = (j_p)_{p \in M}$, and
$k = (k_p)_{p \in N}$.
Assuming $P = M^\uparrow \cup N^\uparrow \cup N \cup N^\downarrow \cup L^\downarrow$, the algebra $\mathcal{T}_{i,j,k}(x)$ is given by
\[
\mathcal{T}_{i,j,k}(x) = 
\bigotimes_{p \in D} J_p
\otimes \bigotimes_{p \in M} \mathcal{A}_{p,j_p}
\otimes \bigotimes_{p \in N} \mathcal{T}_{p,k_p}
\otimes \bigotimes_{p \in L} E_{p,i_p}^*
\otimes \bigotimes_{p \in M^\uparrow \cap L^\downarrow} I_p
\otimes \bigotimes_{p \in U} E_{p,0}^*,
\]
where
$D = N^\downarrow \cup (L^\downarrow \setminus (M \cup M^\uparrow))$
and
$U = N^\uparrow \cup (M^\uparrow \setminus (L \cup L^\downarrow))$.
\end{lem}
\begin{proof}
By Lemma \ref{lem:bij},
there exists a unique generator in \eqref{Tijk}, given by
\[
\mathcal{T}_{i,j,k}'(x) = 
\bigotimes_{p \in D} J_p
\otimes \bigotimes_{p \in M} \mathcal{A}_{p,j_p}
\otimes \bigotimes_{p \in N} \mathcal{T}_{p,k_p}
\otimes \bigotimes_{p \in L} E_{p,i_p}^*
\otimes \bigotimes_{p \in M^\uparrow \cap L^\downarrow} I_p
\otimes \bigotimes_{p \in U} E_{p,0}^*,
\]
where
$D = N^\downarrow \cup (L^\downarrow \setminus (M \cup M^\uparrow))$
and
$U = N^\uparrow \cup (M^\uparrow \setminus (L \cup L^\downarrow))$.
It suffices to demonstrate that $\mathcal{T}_{i,j,k}'(x)$ is closed under matrix multiplication.
This holds true as $\langle J_p \rangle$, $\mathcal{A}_{p,j_p}$ are two-sided ideals in $\mathcal{A}_p$,
$\mathcal{T}_{p,k_p}$ is a two-sided ideal in $\mathcal{T}_p$,
$\langle E_{p,0}^* \rangle$, $\langle E_{p,i_p}^* \rangle$ are two-sided ideals in $\mathcal{A}_p^*$,
and $I_p$ is the identity.
\end{proof}

\begin{prop}\label{prop:Tijk}
With above notation, let $(i,j,k) \in \mathbb{T}$ with
$i = (i_p)_{p \in L}$,
$j = (j_p)_{p \in M}$, and
$k = (k_p)_{p \in N}$.
Define $Q = P \setminus (M^\uparrow \cup N^\uparrow \cup N \cup N^\downarrow \cup L^\downarrow)$.
Assume $Q \neq \emptyset$.
Let $\mathfrak{X}_Q$ denote the generalized wreath product of the association schemes $\{\mathfrak{X}_p\}_{p \in Q}$ over $Q$
and $\mathcal{T}_Q = \mathcal{T}_Q(x_Q)$ denote the Terwilliger algebra of $\mathfrak{X}_Q$ with respect to $x_Q = (x_p)_{p \in Q}$.
Furthermore, let $\mathcal{J}_Q = \mathcal{J}_Q(x_Q)$ denote the principal two-sided ideal of $\mathcal{T}_Q$ generated by the all-ones matrix $J_Q = \bigotimes_{p \in Q}J_p$
Then the algebra $\mathcal{T}_{i,j,k}(x)$ is given by
\[
\mathcal{T}_{i,j,k}(x) = 
\mathcal{J}_Q
\otimes \bigotimes_{p \in D} J_p
\otimes \bigotimes_{p \in M} \mathcal{A}_{p,j_p}
\otimes \bigotimes_{p \in N} \mathcal{T}_{p,k_p}
\otimes \bigotimes_{p \in L} E_{p,i_p}^*
\otimes \bigotimes_{p \in M^\uparrow \cap L^\downarrow} I_p
\otimes \bigotimes_{p \in U} E_{p,0}^*,
\]
where
$D = N^\downarrow \cup (L^\downarrow \setminus (M \cup M^\uparrow))$
and
$U = N^\uparrow \cup (M^\uparrow \setminus (L \cup L^\downarrow))$.
\end{prop}
\begin{proof}
Define
\[
\mathcal{T}_{i,j,k}'(x) = 
\bigotimes_{p \in D} J_p
\otimes \bigotimes_{p \in M} \mathcal{A}_{p,j_p}
\otimes \bigotimes_{p \in N} \mathcal{T}_{p,k_p}
\otimes \bigotimes_{p \in L} E_{p,i_p}^*
\otimes \bigotimes_{p \in M^\uparrow \cap L^\downarrow} I_p
\otimes \bigotimes_{p \in U} E_{p,0}^*,
\]
where
$D = N^\downarrow \cup (L^\downarrow \setminus (M \cup M^\uparrow))$
and
$U = N^\uparrow \cup (M^\uparrow \setminus (L \cup L^\downarrow))$.
According to Lemma \ref{lem:bij},
the generators in \eqref{Tijk} are given by
\[
\left( \bigotimes_{p \in D} J_p \otimes \bigotimes_{p \in H} \mathcal{J}_p \otimes \bigotimes_{p \in U} E_{p,0}^* \right) \otimes \mathcal{T}_{i,j,k}'(x)
\]
for $(D,H,U)$ in the set Lemma \ref{lem:bij} (ii) with respect to $(M,N,L)$.
By Theorem \ref{thm:J} and Lemma \ref{lem:Q=empty}, these generators suffice to produce the algebra
\[
\mathcal{J}_Q \otimes \mathcal{T}_{i,j,k}'(x)
\]
as required.
\end{proof}

\begin{lem}\label{lem:TijkinT}
With above notation, for $(i,j,k) \in \mathbb{T}$, we have $\mathcal{T}_{i,j,k}(x) \subset \mathcal{T}(x)$.
\end{lem}
\begin{proof}
Write
$i = (i_p)_{p \in L}$,
$j = (j_p)_{p \in M}$,
$k = (k_p)_{p \in N}$.
Fix a ordered triple $(D,H,U)$ on the set Lemma \ref{lem:bij} (i) with respect to $(M,N,L)$.
It suffices to show the generator corresponding to $(D,H,U)$ in \eqref{Tijk} contains in $\mathcal{T}(x)$.
By definition, $(D,M,N \cup H,L,M^\uparrow \cap L^\downarrow,U) \in \Omega$.
Since $\mathcal{A}_{p,j_p} \subset \mathcal{A}_p$, $E_{p,i_p}^* \subset \mathcal{A}_p^*$, $\mathcal{T}_{p,k_p} \subset \mathcal{T}_p$, $\mathcal{J}_p \subset \mathcal{T}_p$,
By Lemma \ref{lem:S}, the result follows.
\end{proof}

\begin{lem}\label{lem:idealT}
With above notation, for $(i,j,k) \in \mathbb{T}$, the subspace $\mathcal{T}_{i,j,k}(x)$ constitutes a two-sided ideal of $\mathcal{T}(x)$.
\end{lem}

\begin{proof}
We first demonstrate $\mathcal{T}_{i,j,k}(x)$ is a left ideal of $\mathcal{T}(x)$.
Let
$i = (i_p)_{p \in L}$,
$j = (j_p)_{p \in M}$,
$k = (k_p)_{p \in N}$,
and define
$D = N^\downarrow \cup (L^\downarrow \setminus (M \cup M^\uparrow))$,
$C = M^\uparrow \cap L^\downarrow$,
$U = N^\uparrow \cup (M^\uparrow \setminus (L \cup L^\downarrow))$
and
$Q = P \setminus (M^\uparrow \cup N^\uparrow \cup N \cup N^\downarrow \cup L^\downarrow)$.

Assume $Q = \emptyset$.
According to Proposition \ref{prop:T}, it suffices to show $\mathcal{S}_{\omega}(x)\mathcal{T}_{i,j,k}(x) \subset \mathcal{T}_{i,j,k}(x)$ for $\omega = (D',M',N',L',C',U') \in \Omega$,
where
\[
\mathcal{S}_{\omega}(x)= 
\bigotimes_{p \in D'} J_p
\otimes \bigotimes_{p \in M'} \mathcal{A}_p
\otimes \bigotimes_{p \in N'} \mathcal{T}_p
\otimes \bigotimes_{p \in L'} \mathcal{A}_p^*
\otimes \bigotimes_{p \in C'} I_p
\otimes \bigotimes_{p \in U'} E_{p,0}^*.
\]
If there exists $p \in M \cap D'$, since $j_p \neq 0$, 
the $p$-th factor of $\mathcal{S}_{\omega}(x)\mathcal{T}_{i,j,k}(x)$ is $J_p \mathcal{A}_{p,j_p} = 0$.
If there exists $p \in N \cap D'$, since $k_p \neq 0$, 
the $p$-th factor of $\mathcal{S}_{\omega}(x)\mathcal{T}_{i,j,k}(x)$ is $J_p \mathcal{T}_{p,k_p} = 0$.
If there exists $p \in (L \cup C \cup U) \cap (D' \cup M' \cup N')$, since $L \cup C \cup U \subset M^\uparrow \cup N^\uparrow$, there exists $q \in M \cup N$ such that $q < p$.
Since $(D', M' \cup N', L' \cup C' \cup U') \in \Delta$, this $q$ satisfies $q \in D'$.
Thus, there exists $q \in (M \cup N) \cap D'$ in this case.
Therefore, we have shown that $\mathcal{S}_{\omega}(x)\mathcal{T}_{i,j,k}(x) = 0 \subset \mathcal{T}_{i,j,k}(x)$ if $(M \cup N) \cap D' \neq \emptyset$ or $(L \cup C \cup U) \cap (D' \cup M' \cup N') \neq \emptyset$.
Similarly, one can demonstrate that $\mathcal{S}_{\omega}(x)\mathcal{T}_{i,j,k}(x) = 0 \subset \mathcal{T}_{i,j,k}(x)$ if $(N \cup L) \cap U' \neq \emptyset$ or $(M \cup C \cup D) \cap (N' \cup L' \cup U') \neq \emptyset$.
Consider the case where $(M \cup N) \cap D' = (L \cup C \cup U) \cap (D' \cup M' \cup N') = (N \cup L) \cap U' = (M \cup C \cup D) \cap (N' \cup L' \cup U') = \emptyset$.
If $p \in D \cup M$, we have $p \in D' \cup M' \cup C'$ and so
the $p$-th factor of $\mathcal{S}_{\omega}(x)$ is contained in $\mathcal{A}_p$, while that of $\mathcal{T}_{i,j,k}(x)$ is $J_p$ or $\mathcal{A}_{p,j_p}$, which are two-sided ideals of $\mathcal{A}_p$.
If $p \in N$, 
the $p$-th factor of $\mathcal{T}_{i,j,k}(x)$ is $\mathcal{T}_{p,k_p}$, which is a two-sided ideal of $\mathcal{T}_p$.
If $p \in L \cup U$, we have $p \in L' \cup C' \cup U'$ and so
the $p$-th factor of $\mathcal{S}_{\omega}(x)$ is contained in $\mathcal{A}_p^*$, while that of $\mathcal{T}_{i,j,k}(x)$ is $E_{p,0}^*$ or $E_{p,i_p}^*$, which are two-sided ideals of $\mathcal{A}_p^*$.
If $p \in C$, we have $p \in C'$ and so 
the $p$-th factors of $\mathcal{S}_{\omega}(x)$ and $\mathcal{T}_{i,j,k}(x)$ are the identity $I_p$.
Therefore we conclude that $\mathcal{S}_{\omega}(x)\mathcal{T}_{i,j,k}(x) \subset \mathcal{T}_{i,j,k}(x)$ in all cases.

Assume $Q \neq \emptyset$ and set $Q' = P \setminus Q$.
Define
\[
\mathcal{T}_{i,j,k}'(x) = 
\bigotimes_{p \in D} J_p
\otimes \bigotimes_{p \in M} \mathcal{A}_{p,j_p}
\otimes \bigotimes_{p \in N} \mathcal{T}_{p,k_p}
\otimes \bigotimes_{p \in L} E_{p,i_p}^*
\otimes \bigotimes_{p \in M^\uparrow \cap L^\downarrow} I_p
\otimes \bigotimes_{p \in U} E_{p,0}^*,
\]
where
$D = N^\downarrow \cup (L^\downarrow \setminus (M \cup M^\uparrow))$
and
$U = N^\uparrow \cup (M^\uparrow \setminus (L \cup L^\downarrow))$.
By Proposition \ref{prop:Tijk}, we have $\mathcal{T}_{i,j,k}(x) = \mathcal{J}_Q \otimes \mathcal{T}_{i,j,k}'(x)$.
By definition, $\mathcal{J}_Q$ is a two-sided ideal of $\mathcal{T}_Q$
and as we have just demonstrated, $\mathcal{T}_{i,j,k}'(x)$ is a left ideal of $\mathcal{T}_{Q'}$.
Therefore, by Corollary \ref{cor:T}, $\mathcal{T}_{i,j,k}(x)$ is also a left ideal of $\mathcal{T}(x)$.

Thus, we have shown that $\mathcal{T}_{i,j,k}(x)$ is a left ideal of $\mathcal{T}(x)$ in both scenarios.
Similarly, one can establish that $\mathcal{T}_{i,j,k}(x)$ is a right ideal of $\mathcal{T}(x)$, and hence the result follows.
\end{proof}

\begin{lem}\label{lem:orthogonalT}
With above notation, the subspaces $\{\mathcal{T}_{i,j,k}(x)\}_{(i,j,k) \in \mathbb{T}}$ are mutually orthogonal.
\end{lem}
\begin{proof}
For $(i,j,k)\in \mathbb{T}$ with $i = (i_p)_{p \in L}$,
$j = (j_p)_{p \in M}$,
$k = (k_p)_{p \in N}$,
and for $(D,H,U)$ in the set of Lemma \ref{lem:bij} (i) with respect to $(M,N,L)$, we define
\[
\mathcal{T}_{i,j,k,D,H,U}'(x) = 
\bigotimes_{p \in D} J_p
\otimes \bigotimes_{p \in M} \mathcal{A}_{p,j_p}
\otimes \bigotimes_{p \in N} \mathcal{T}_{p,k_p}
\otimes \bigotimes_{p \in H} \mathcal{J}_p
\otimes \bigotimes_{p \in L} E_{p,i_p}^*
\otimes \bigotimes_{p \in M^\uparrow \cap L^\downarrow} I_p
\otimes \bigotimes_{p \in U} E_{p,0}^*.
\]
Let $(i,j,k)\in \mathbb{T}$ with $i = (i_p)_{p \in L}$,
$j = (j_p)_{p \in M}$,
$k = (k_p)_{p \in N}$
and let
$(\tilde{i}, \tilde{j}, \tilde{k}) \in \mathbb{T}$ with $\tilde{i} = (\tilde{i}_p)_{p \in \tilde{L}}$,
$\tilde{j} = (\tilde{j}_p)_{p \in \tilde{M}}$,
$\tilde{k} = (\tilde{k}_p)_{p \in \tilde{N}}$.
Assume $(i,j,k) \neq (\tilde{i}, \tilde{j}, \tilde{k})$.

Take $(D,H,U)$, $(\tilde{D},\tilde{H},\tilde{U})$ in the set Lemma \ref{lem:bij} (i) with respect to $(M,N,L)$, $(\tilde{M},\tilde{N},\tilde{L})$, respectively. 
We aim to show $\mathcal{T}_{i,j,k,D,H,U}'(x)\mathcal{T}_{\tilde{i},\tilde{j},\tilde{k},\tilde{D},\tilde{H},\tilde{U}}'(x) = 0$.
We abbreviate $\mathcal{T} = \mathcal{T}_{i,j,k,D,H,U}'(x)$ and $\tilde{\mathcal{T}} = \mathcal{T}_{\tilde{i},\tilde{j},\tilde{k},\tilde{D},\tilde{H},\tilde{U}}'(x)$.

If there exists $p \in M \cap \tilde{D}$,
then the $p$-th factor of $\mathcal{T}\tilde{\mathcal{T}}$
is $\mathcal{A}_{p,j_p}J_p = 0$ since $j_p \neq 0$.
Similarly, if there exists $p \in D \cap \tilde{M}$,
the $p$-th factor of $\mathcal{T}\tilde{\mathcal{T}}$ is zero.
If there exists $p \in N \cap \tilde{D}$,
then the $p$-th factor of $\mathcal{T}\tilde{\mathcal{T}}$
is $\mathcal{T}_{p,k_p}J_p = 0$ since $k_p \neq 0$.
Similarly, if there exists $p \in D \cap \tilde{N}$,
the $p$-th factor of $\mathcal{T}\tilde{\mathcal{T}}$ is zero.
If there exists $p \in N \cap \tilde{H}$,
then the $p$-th factor of $\mathcal{T}\tilde{\mathcal{T}}$
is $\mathcal{T}_{p,k_p}\mathcal{J}_p = 0$ since $k_p \neq 0$.
Similarly, if there exists $p \in H \cap \tilde{N}$,
the $p$-th factor of $\mathcal{T}\tilde{\mathcal{T}}$ is zero.
If there exists $p \in N \cap \tilde{U}$,
then the $p$-th factor of $\mathcal{T}\tilde{\mathcal{T}}$
is $\mathcal{T}_{p,k_p}E_{p,0}^* = 0$ since $k_p \neq 0$.
Similarly, if there exists $p \in U \cap \tilde{N}$,
the $p$-th factor of $\mathcal{T}\tilde{\mathcal{T}}$ is zero.
If there exists $p \in L \cap \tilde{U}$,
then the $p$-th factor of $\mathcal{T}\tilde{\mathcal{T}}$
is $E_{p,i_p}^*E_{p,0}^* = 0$ since $i_p \neq 0$.
Similarly, if there exists $p \in U \cap \tilde{L}$,
the $p$-th factor of $\mathcal{T}\tilde{\mathcal{T}}$ is zero.
Therefore, it suffices to show at least one of the following sets is non-empty:
$M \cap \tilde{D}$, $D \cap \tilde{M}$, $N \cap \tilde{D}$, $D \cap \tilde{N}$, $N \cap \tilde{H}$, $H \cap \tilde{N}$,
$N \cap \tilde{U}$, $U \cap \tilde{N}$, $L \cap \tilde{U}$, $U \cap \tilde{L}$.

Assume $M \not\subset \tilde{M}$.
In this case, there exists $p \in M$ such that $p \not\in \tilde{M}$.
Then $p$ is either in $\tilde{D}$, $\tilde{M}^\uparrow \cap \tilde{L}^\downarrow$, or $\tilde{N} \cup \tilde{H} \cup \tilde{L} \cup \tilde{U}$.
If $p \in \tilde{D}$, then $p \in M \cap \tilde{D}$.
If $p \in \tilde{M}^\uparrow \cap \tilde{L}^\downarrow$, then there exists $q \in D \cap \tilde{M}$ such that $q < p$.
If $p \in \tilde{N} \cup \tilde{H} \cup \tilde{L} \cup \tilde{U}$, then there exists $q \in L \cap \tilde{U}$ such that $p < q$.
Therefore, the assertion holds if $M \not\subset \tilde{M}$.
Similarly, the assertion holds if $\tilde{M} \not \subset M$.

Assume $N \not\subset \tilde{N}$.
In this case, there exists $p \in N$ such that $p \not\in \tilde{N}$.
Then $p$ is either in $\tilde{D}$, $\tilde{M}$, $\tilde{H}$, $\tilde{L}$, $\tilde{M}^\uparrow \cap \tilde{L}^\downarrow$, or $\tilde{U}$.
If $p \in \tilde{D}$, then $p \in N \cap \tilde{D}$.
If $p \in \tilde{M}$, then there exists $q \in U \cap \tilde{L}$ such that $p < q$.
If $p \in \tilde{H}$, then $p \in N \cap \tilde{H}$.
If $p \in \tilde{L}$, then there exists $q \in D \cap \tilde{M}$ such that $q < p$.
If $p \in \tilde{M}^\uparrow \cap \tilde{L}^\downarrow$, then there exists $q \in D \cap \tilde{M}$ such that $q < p$.
If $p \in \tilde{U}$, then $p \in N \cap \tilde{U}$.
Therefore, the assertion holds if $N \not\subset \tilde{N}$.
Similarly, the assertion holds if $\tilde{N} \not \subset N$.

Assume $L \not\subset \tilde{L}$.
In this case, there exists $p \in L$ such that $p \not\in \tilde{L}$.
Then $p$ is either $\tilde{D} \cup \tilde{M} \cup \tilde{N} \cup \tilde{H}$, $\tilde{M}^\uparrow \cap \tilde{L}^\downarrow$, or $\tilde{U}$.
If $p \in \tilde\tilde{D} \cup \tilde{M} \cup \tilde{N} \cup \tilde{H}$, then there exists $q \in M \cap \tilde{D}$ such that $q < p$.
If $p \in \tilde{M}^\uparrow \cap \tilde{L}^\downarrow$, then there exists $q \in U \cap \tilde{L}$ such that $p < q$.
If $p \in \tilde{U}$, then $p \in L \cap \tilde{U}$.
Therefore, the assertion holds if $L \not\subset \tilde{L}$.
Similarly, the assertion holds if $\tilde{L} \not \subset L$.

Consider the case $(M,N,L) = (\tilde{M},\tilde{N},\tilde{L})$.
Since $(i,j,k) \neq (\tilde{i}, \tilde{j}, \tilde{k})$,
there exists either
a $p \in M = \tilde{M}$ such that $j_p \neq \tilde{j}_p$,
a $p \in N = \tilde{N}$ such that $i_p \neq \tilde{i}_p$,
or 
a $p \in L = \tilde{L}$ such that $k_p \neq \tilde{k}_p$.
If there exists $p \in M = \tilde{M}$ such that $j_p \neq \tilde{j}_p$,
the $p$-th factor of $\mathcal{T}\tilde{\mathcal{T}}$
is $\mathcal{A}_{p,j_p}\mathcal{A}_{p,\tilde{j}_p} = 0$.
Similarly, the assertion holds for the other two cases.
\end{proof}

\begin{lem}\label{lem:AinT}
With above notation, we have $\mathcal{A} \subset \sum_{(i,j,k) \in \mathbb{T}}\mathcal{T}_{i,j,k}(x)$.
\end{lem}
\begin{proof}
By Theorem \ref{thm:A}, it suffices to demonstrate that $\mathcal{A}_j \subset \sum_{(i,j,k) \in \mathbb{T}}\mathcal{T}_{i,j,k}(x)$ for $j \in \mathbb{J}$.
If $j = 0$,
we have $\mathcal{A}_0 = \langle J_X \rangle \subset \mathcal{J}(x) = \mathcal{T}_{0,0,0}$.
Assume $j \neq 0$.
Let $j = (j_p)_{p \in M} \in \mathbb{J}$ and let $B_{p,j_p} \in \mathcal{A}_{p,j_p}$.
We need to verify whether
\[
\bigotimes_{p \not\in M \cup M^\uparrow} J_p \otimes \bigotimes_{p \in M} B_{p,j_p} \otimes \bigotimes_{p \in M^\uparrow} I_p
\in
\sum_{(i,j,k) \in \mathbb{T}}\mathcal{T}_{i,j,k}(x).
\]
By the decomposition \eqref{Tp}, without loss of generality, we may assume there exists $0 \le k_p \le t_p$ such that $B_{p,j_p} \in \mathcal{T}_{p,k_p}$ for $p \in M$.
Define
$N = \{p \in M \mid k_p \ge 1\}$ and 
$H = \{p \in M \mid k_p = 0\}$.
Set $\mathbb{I}(M^\uparrow) = \{ (i_p)_{p \in L} \in \mathbb{I} \mid L \subset M^\uparrow\}$.
For each $i = (i_p)_{p \in L} \in \mathbb{I}(M^\uparrow)$, define $E_i^*(M^\uparrow)$ analogously to $E_i^*(x)$ as in \eqref{E*} but with indices restricted to $M^\uparrow$.
More precisely,
\[
E_i^*(M^\uparrow) = \bigotimes_{p \in M^\uparrow \cap L^\downarrow} I_p  \otimes \bigotimes_{p \in L} E_{p,i_p}^* \otimes \bigotimes_{p \in M^\uparrow \setminus (L \cup L^\downarrow)} E_{p,0}^*.
\]
Then $\{E_i^*(M^\uparrow)\}_{i \in \mathbb{I}(M^\uparrow)}$ constitutes the complete set of dual primitive idempotents of the generalized wreath product of $\{\mathfrak{X}_p\}_{p \in M^\uparrow}$ over $M^\uparrow$ with respect to $(x_p)_{p \in M^\uparrow}$.
In particular, the sum of all $E_i^*(M^\uparrow)$ is the identity.
For $i = (i_p)_{p \in L} \in \mathbb{I}(M^\uparrow)$, we have $(i,j,k) \in \mathbb{T}$ where $j = (j_p)_{p \in M \cap L^\downarrow}$,
$k = (k_p)_{p \in N \setminus L^\downarrow}$.
Additionally, $(D,H',U)$ is in the set of Lemma \ref{lem:bij} (i) with respect to $(M \cap L^\downarrow, N \setminus L^\downarrow, L)$, where $D = P \setminus (M \cup M^\uparrow)$, $H' = H \setminus L^\downarrow$, and $U = M^\uparrow \setminus (L \cup L^\downarrow)$.
Thus
\begin{align*}
&\bigotimes_{p \not\in M \cup M^\uparrow} J_p \otimes \bigotimes_{p \in M} B_{p,j_p} \otimes \bigotimes_{p \in M^\uparrow} I_p \\
&=
\bigotimes_{p \not\in M \cup M^\uparrow} J_p \otimes \bigotimes_{p \in M} B_{p,j_p} \otimes \left( \sum_{i \in \mathbb{I}(M^\uparrow)} E_i^*(M^\uparrow) \right) \\
&= \sum_{i = (i_p)_{p \in L} \in \mathbb{I}(M^\uparrow)}
\bigotimes_{p \not\in M \cup M^\uparrow} J_p
\otimes \bigotimes_{p \in M \cap L^\downarrow} B_{p,j_p} 
\otimes \bigotimes_{p \in N \setminus L^\downarrow} B_{p,j_p} 
\otimes \bigotimes_{p \in H \setminus L^\downarrow} B_{p,j_p} \\
&\phantom{= \sum_{i = (i_p)_{p \in L} \in \mathbb{I}(M^\uparrow)} \bigotimes}
\otimes \bigotimes_{p \in M^\uparrow \cap L^\downarrow}  I_p
\otimes \bigotimes_{p \in L} E_{p,i_p}^*
\otimes \bigotimes_{p \in M^\uparrow \setminus (L \cup L^\downarrow)} E_{p,0}^*.
\end{align*}
Since
$B_{p,j_p} \subset \mathcal{A}_{p,j_p}$ if $p \in M \cap L^\downarrow$,
$B_{p,j_p} \subset \mathcal{T}_{p,k_p}$ if $p \in N \setminus L^\downarrow$, and
$B_{p,j_p} \subset \mathcal{J}_p$ if $p \in H \setminus L^\downarrow$,
the summand corresponding to $i \in \mathbb{I}(M^\uparrow)$ belongs to $\mathcal{T}_{i,j,k}(x)$, where $j = (j_p)_{p \in M \cap L^\downarrow}$, $k = (k_p)_{p \in N \setminus L^\downarrow}$.
The result follows.
\end{proof}

\begin{lem}\label{lem:A*inT}
With above notation, we have $\mathcal{A}^*(x) \subset \sum_{(i,j,k) \in \mathbb{T}}\mathcal{T}_{i,j,k}(x)$.
\end{lem}
\begin{proof}
Analogous to Lemma \ref{lem:AinT}.
\end{proof}

\begin{thm}\label{thm:T}
With above notation,
the Terwilliger algebra $\mathcal{T}(x)$ of $\mathfrak{X}$ with respcet to $x$ is represented as a direct sum of central simple algebras:
\[
\mathcal{T}(x) = \bigoplus_{(i,j,k) \in \mathbb{T}} \mathcal{T}_{i,j,k}(x),
\]
where $\mathcal{T}_{i,j,k}(x)$ is specified in \eqref{Tijk}.
\end{thm}
\begin{proof}
First, we establish that each $\mathcal{T}_{i,j,k}(x)$ is a central simple algebra.
For $(i,j,k) \in \mathbb{T}$ with 
$i = (i_p)_{p \in L}$,
$j = (j_p)_{p \in M}$, and
$k = (k_p)_{p \in N}$,
by Propsition \ref{prop:Tijk},
we have, as algebras,
\[
\mathcal{T}_{i,j,k}(x) \cong \mathcal{J}_Q \otimes \bigotimes_{p \in M} \mathcal{A}_{p,j_p}
\otimes \bigotimes_{p \in N} \mathcal{T}_{p,k_p},
\]
where $Q = P \setminus (M^\uparrow \cup N^\uparrow \cup N \cup N^\downarrow \cup L^\downarrow)$.
Thus, $\mathcal{T}_{i,j,k}(x)$ is a central simple algebra since it is a tensor product of central simple algebras.
Next, we confirm the desired equality.
By Lemmas \ref{lem:TijkinT}, \ref{lem:AinT} and \ref{lem:A*inT}, it follows that
\[
\mathcal{A}, \mathcal{A}^*(x) \subset \sum_{(i,j,k) \in \mathbb{T}} \mathcal{T}_{i,j,k}(x) \subset \mathcal{T}(x).
\]
According to Lemma \ref{lem:orthogonalT}, the sum is closed under matrix multiplication, this proves that
the sum coincides with $\mathcal{T}(x)$, given that $\mathcal{T}(x)$ is generated by $\mathcal{A}$ and $\mathcal{A}^*(x)$.
Finally, in accordance with Lemmas \ref{lem:idealT} and \ref{lem:orthogonalT}, the sum is direct.
\end{proof}

\begin{thm}\label{thm:epsilon}
For each $p \in P$, let $\{\varepsilon_{p,k}(x_p)\}_{k=0}^{t_p}$ denote the primitive central idempotents of the Terwilliger algebra of $\mathfrak{X}_p$
with $\varepsilon_{p,0}(x_p) = \epsilon_p$, the principal idempotent.
We abbreviate $\varepsilon_{p,k} = \varepsilon_{p,k}(x_p)$.
Then the primitive central idempotents of the Terwilliger algebra of $\mathfrak{X}$ are given by
\begin{align*}
\epsilon_Q
\otimes \bigotimes_{p \in D} E_{p,0}
\otimes \bigotimes_{p \in M} E_{p,j_p}
\otimes \bigotimes_{p \in N} \varepsilon_{p,k_p}
\otimes \bigotimes_{p \in L} E_{p,i_p}^*
\otimes \bigotimes_{p \in M^\uparrow \cap L^\downarrow} I_p
\otimes \bigotimes_{p \in U} E_{p,0}^*,
&&
(i,j,k) \in \mathbb{T},
\end{align*}
where 
$i = (i_p)_{p \in L}$,
$j = (j_p)_{p \in M}$,
$k = (k_p)_{p \in N}$,
$D = N^\downarrow \cup (L^\downarrow \setminus (M \cup M^\uparrow))$,
$U = N^\uparrow \cup (M^\uparrow \setminus (L \cup L^\downarrow))$,
$Q = P \setminus (M^\uparrow \cup N^\uparrow \cup N \cup N^\downarrow \cup L^\downarrow)$.
Here, $\mathfrak{X}_Q$ denotes the generalized wreath product of the association schemes $\{\mathfrak{X}_p\}_{p \in Q}$ over $Q$
and $\mathcal{T}_Q(x_Q)$ denotes the Terwilliger algebra of $\mathfrak{X}_Q$ with respect to $x_Q = (x_p)_{p \in Q}$
 with $\epsilon_Q$ representing the principal idempotent of $\mathcal{T}_Q(x_Q)$.
\end{thm}

\begin{proof}
By Theorem \ref{thm:T}, every primitive central idempotent is the identity in $\mathcal{T}_{i,j,k}(x)$ for some $(i,j,k) \in \mathbb{T}$.
The result follows from Proposition \ref{prop:Tijk} and the definitions of $\{E_{p,j}\}_{j=0}^{l_p}$, $\{\varepsilon_{p,k}(x_p)\}_{k=0}^{t_p}$ and $\epsilon_Q$.
\end{proof}

\section{Special cases}
In this section, we discuss applications of Theorems \ref{thm:J}, \ref{thm:eps}, \ref{thm:T} and \ref{thm:epsilon} to the special cases studied in \cite{HKM, MX}.
We use the following notation.
Let $\mathfrak{X} = (X, \{R_i\}_{i=0}^d)$ and $\mathfrak{Y} = (Y, \{S_i\}_{i=0}^e)$ denote association schemes.
Let $\mathcal{A}$ (resp.\ $\mathcal{B}$) denote the adjacency algebra of $\mathfrak{X}$ (resp.\ $\mathfrak{Y}$).
The direct sum decompositions into central simple algebras are given by
\begin{align*}
\mathcal{A} = \bigoplus_{i=0}^l \mathcal{A}_i, && \mathcal{B} = \bigoplus_{j=0}^k \mathcal{B}_j.
\end{align*}
As usual, we assume $\mathcal{A}_0$ (resp.\ $\mathcal{B}_0$) is the principal two-sided ideal of $\mathcal{A}$ (resp.\ $\mathcal{B}$) spanned by the all-ones matrix.
Let $\{E_i(x)\}_{i=0}^l$ (resp.\ $\{F_j(y)\}_{j=0}^k$) denote the primitive central idempotent of $\mathcal{A}$ (resp.\ $\mathcal{B}$).
As usual, we assume $E_0 = |X|^{-1}J_X$ and $F_0 = |Y|^{-1}J_Y$.
Fix $x \in X$ and $y \in Y$.
Let $\mathcal{A}^*(x)$ (resp.\ $\mathcal{B}^*(y)$) denote the dual adjacency algebra of $\mathfrak{X}$ (resp.\ $\mathfrak{Y}$) with respect to $x$ (resp.\ $y$) 
and let $\{E_i^*(x)\}_{i=0}^d$ (resp.\ $\{F_j^*(y)\}_{j=0}^e$) denote the primitive idempotents of $\mathfrak{X}$ (resp.\ $\mathfrak{Y}$) with respect to $x$ (resp.\ $y$).
Let $\mathcal{T}(x)$ (resp.\ $\mathcal{S}(y)$) denote the Terwilliger algebra of $\mathfrak{X}$ (resp.\ $\mathfrak{Y}$) with respect to $x$ (resp.\ $y$).
The direct sum decompositions into central simple algebras are given by
\begin{align*}
\mathcal{T}(x) = \bigoplus_{i=0}^t \mathcal{T}_i(x), && \mathcal{S}(y) = \bigoplus_{j=0}^s \mathcal{S}_j(y).
\end{align*}
As usual, we assume $\mathcal{T}_0(x)$ (resp.\ $\mathcal{S}_0(y)$) is the principal two-sided ideal of $\mathcal{T}(x)$ (resp.\ $\mathcal{S}(y)$) generated by the all-ones matrix.
Let $\{\varepsilon_i(x)\}_{i=0}^t$ (resp.\ $\{\varrho_j(y)\}_{j=0}^s$) denote the primitive central idempotent of $\mathcal{T}(x)$ (resp.\ $\mathcal{S}(y)$).
As usual, we assume $\varepsilon_0(x)$ and $\varrho_0(y)$ are the principal idempotents.

\subsection{Direct products}

Let $P$ be a poset consisting of two incomparable elements.
Then the generalized wreath product over $P$ is called the direct product of two association schemes.
The direct product of two association schemes $\mathfrak{X}$ and $\mathfrak{Y}$ is denoted by $\mathfrak{X} \times \mathfrak{Y}$.

\begin{cor}[{\cite[Theorem 3.1]{HKM}}]
With above notation, the Terwilliger algebra of $\mathfrak{X} \times \mathfrak{Y}$ with respect to $(x,y)$ is represented as a direct sum of central simple algebras:
\[
\bigoplus_{i=0}^t \bigoplus_{j=0}^s \left( \mathcal{T}_i(x) \otimes \mathcal{S}_j(y) \right).
\]
Moreover, $\mathcal{T}_0(x) \otimes \mathcal{S}_0(y)$ is the principal two-sided ideal generated by the all-ones matrix.
\end{cor}

\begin{cor}[c.f.\ {\cite[Theorem 3.1]{HKM}}]
With above notation, the primitive central idempotents of the Terwilliger algebra of $\mathfrak{X} \times \mathfrak{Y}$ with respect to $(x,y)$ are given by
\begin{align*}
\varepsilon_i(x) \otimes \varrho_j(y),
&&
0 \le i \le t, \; 0 \le j \le s.
\end{align*}
Moreover, $\varepsilon_0(x) \otimes \varrho_0(y)$ is the principal idempotent.
\end{cor}

\subsection{Wreath products}

Let $P$ be a poset consisting of two comparable elements.
Then the generalized wreath product over $P$ is called the wreath product of two association schemes.
The wreath product of two association schemes $\mathfrak{X}$ and $\mathfrak{Y}$ is denoted by $\mathfrak{X} \wr \mathfrak{Y}$
if, in the poset $P$, the corresponding element of $\mathfrak{X}$ is less than that of $\mathfrak{Y}$.

\begin{cor}[{\cite[Lemma 5.1 (i)]{MX}}]
With above notation, the principal two-sided ideal of the Terwilliger algebra of $\mathfrak{X} \wr \mathfrak{Y}$ with respect to $(x,y)$ generated by the all-ones matrix is given by
\[
J_X \otimes \mathcal{S}_0(x) + 
\mathcal{A}^*(x)J_X \otimes F_0^*(y)\mathcal{B} +
J_X\mathcal{A}^*(x) \otimes \mathcal{B}F_0^*(y) +
\mathcal{T}_0(x) \otimes F_0^*(y).
\]
\end{cor}

\begin{cor}[{\cite[Proposition 5.2]{MX}}]
With above notation, the Terwilliger algebra of $\mathfrak{X} \wr \mathfrak{Y}$ with respect to $(x,y)$ is represented as a direct sum of central simple algebras:
\[
\mathcal{J}(x,y) \oplus \bigoplus_{j=1}^s \left( J_X \otimes \mathcal{S}_j(y)\right) \oplus \bigoplus_{i=1}^t \left( \mathcal{T}_i(x) \otimes F_0^*(y)\right) \oplus \bigoplus_{i=1}^d\bigoplus_{j=1}^e \left(\mathcal{A}_i \otimes F_j^*(y)\right),
\]
where $\mathcal{J}(x,y)$ is the principal two-sided ideals generated by the all-ones matrix.
\end{cor}
\begin{cor}[{\cite[Theorem 5.3]{MX}}]
With above notation, the primitive central idempotents of the Terwilliger algebra of $\mathfrak{X} \wr \mathfrak{Y}$ are given by
\[
\epsilon(x,y) = E_0^*(x) \otimes F_0^*(y) + E_0 \otimes \left( \varrho_0(y) - F_0^*(y)\right) + \left( \varepsilon_0(x) - E_0^*(x)\right) \otimes F_0^*(y),
\]
\begin{align*}
E_0 \otimes \varrho_j(y) && 1 \le j \le s,
\end{align*}
\begin{align*}
\varepsilon_i(x) \otimes F_0^*(y) && 1 \le i \le t,
\end{align*}
\begin{align*}
E_i \otimes F_j^*(y) && 1 \le i \le l, \; 1 \le j \le e.
\end{align*}
where $\epsilon(x,y)$ is the principal idempotent.

\end{cor}

\section*{Acknowledgement}
This work was supported by JSPS KAKENHI Grant Number JP23K12953.

\bigskip

\noindent
Yuta Watanabe \\
Department of Mathematics Education \\
Aichi University of Education \\
1 Hirosawa, Igaya-cho, Kariya, Aichi 448-8542, Japan. \\
email: \texttt{ywatanabe@auecc.aichi-edu.ac.jp}
\end{document}